\documentclass{amsart}
\usepackage{amsmath}
\usepackage{amsfonts}
\usepackage{amssymb}
\usepackage{graphicx}
\usepackage{psfrag}

\newtheorem{theorem}{Theorem}
\newtheorem{lemma}[theorem]{Lemma}
\newtheorem{proposition}[theorem]{Proposition}

\newcommand{\pd}    {\partial}
\newcommand{\fh}    {{\mathcal{H}}}

\newcommand{\cL}    {\mathcal{L}}

\newcommand{\cO}    {\mathcal{O}}
\newcommand{\cP}    {\mathcal{P}}
\newcommand{\cS}    {\mathcal{S}}
\newcommand{\cA}    {\mathcal{A}}
\newcommand{\R}     {\mathbb{R}}

\renewcommand{\l}   {\left}
\renewcommand{\r}   {\right}

\newcommand{\hq}{{\hat{Q}}}

\newcommand{\llangle}{\langle\!\langle}
\newcommand{\rrangle}{\rangle\!\rangle}

\begin{document}

\title[Asymptotics of the Ricci flow neckpinch]{Precise asymptotics of the Ricci flow neckpinch}

\address[Sigurd Angenent]{ University of Wisconsin - Madison}

\email{angenent@math.wisc.edu}

\urladdr{http://www.math.wisc.edu/\symbol{126}angenent/}

\author{Sigurd Angenent \& Dan Knopf}

  \thanks{First author partially supported by NSF grant DMS-0101124.
   Second author partially supported by NSF grants DMS-0511184, DMS-0505920,
   and a University of Texas Summer Research Assignment.}

  \address[Dan Knopf]{University of Texas at Austin}
  \email{danknopf@math.utexas.edu}
  \urladdr{http://www.ma.utexas.edu/\symbol{126}danknopf} \maketitle

\tableofcontents{}

\setlength{\parskip}{1ex plus 1ex}

\section{Introduction}

\subsection{Antecedents}

In virtually all known applications of the Ricci flow, it is valuable to have a good
understanding of singularity formation. Heuristically, there are at least three
reasons for this. The first is that one expects finite-time singularities to form for
a broad spectrum of initial data. Indeed, such singularities are inevitable if the
scalar curvature is strictly positive. The second reason is that one expects the
geometry of a solution to resemble a standard model (for example, a self-similar
solution) in a space-time neighborhood of a developing singularity. The third reason
is that having a sufficiently detailed picture of a developing singularity
facilitates the geometric-topological surgeries by which Ricci flow decomposes a
given manifold.

Whenever a compact solution $(M^{n},g(\cdot))$ of Ricci flow encounters a singularity
at time $T<\infty$, standard short-time existence results imply that
\[
\lim_{t\nearrow T}\max_{x\in M^{n}}|\operatorname*{Rm}(x,t)|=\infty.
\]
Recent results of N.~Sesum \cite{Se} allow one to replace this condition by
\[
\limsup_{t\nearrow T}\max_{x\in M^{n}}|\operatorname*{Rc}(x,t)|=\infty.
\]
The most interesting cases are those where a \emph{local singularity }forms, that is,
where there exists an open set $\Omega\subset M^{n}$ such that
\[
\sup_{\Omega\times\lbrack0,T)}|\operatorname*{Rc}(x,t)|<\infty.
\]
The first rigorous constructions of local singularities were done by M.~Simon
\cite{Si}. Here the manifold is a noncompact warped product $\mathbb{R}
\times_{f}S^{n}$, and a supersolution of the Ricci flow \textsc{pde} is used to prove
that $f$ vanishes in finite time on a proper subset of $\mathbb{R}$.  A second class
of examples was constructed in \cite{FIK}. Here the manifold is a complex line bundle
$\mathbb{C\hookrightarrow L}_{-k}^{n}\rightarrow \mathbb{CP}^{n-1}$ with $1\leq k\leq
n-1$. As the singularity forms, the flow performs an algebraic-geometric blow-down of
the $\mathbb{CP}^{n-1}$, while the rest of the manifold converges locally smoothly to
a metric cone on $(\mathbb{C}^{n}\backslash\{0\})/\mathbb{Z}_{k}$.

As G.~Perelman writes, \textquotedblleft The most natural way of forming a
singularity in finite time is by pinching an (almost) round cylindrical
neck.\textquotedblright\ \cite{P1} This is the kind of local singularity that we
analyze in this paper. It is a continuation of our earlier work \cite{AK}, where we
gave the first rigorous examples of local singularity formation on compact manifolds
by constructing neckpinches for rotationally-symmetric metrics on $S^{n+1}$. In that
paper, we obtained local \emph{a priori }estimates for the space-time scales at which
a developing neckpinch singularity resembles the self-similarly-shrinking cylinder
soliton. The present paper shows that those estimates are sharp, and gives precise
asymptotics for neckpinch formation.

When considering a local singularity, two natural and important questions arise. (1)
What is the nature of the set of points in space on which the metric becomes
singular?\ (2) What is the asymptotic behavior of the solution near this set as the
singularity time is approached? In the past two decades, a rich literature of both
rigorous and formal matched asymptotics has developed for analyzing the local
behavior of singular solutions of nonlinear \textsc{pde} such as $u_{t}=\Delta\left(
  \log u\right) $ or $u_{t}=\Delta u+F(u)$, where $F(u)=u^{p}$ or $F(u)=e^{u}$. Some
of the many noteworthy results of this type are \cite{Weissler}, \cite{MW},
\cite{FM}, \cite{GP}, \cite{GK1, GK2, GK3}, \cite{HV1, HV2, HV3}, \cite{King}, and
\cite{FK}. A few results specific to geometric evolution equations such as mean
curvature flow and harmonic map flow are \cite{AAG}, \cite{AV0, AV}, \cite{VHK},
\cite{AA}, and \cite{CK}.

\subsection{Overview}

This paper is divided into two main parts.

\medskip

\noindent\textbf{Rigorous asymptotics. }In Section \ref{Round}, we derive
rigorous asymptotics for rotationally symmetric neckpinches. The set-up is as
follows. Let $g_{\operatorname*{can}}$ denote the round metric of radius $1$ on
$S^{n}$. Then any $\operatorname{SO}(n+1)$-invariant metric on $S^{n+1}$ can be
written as
\[
g=\varphi^{2}\,dx^{2}+\psi^{2}\,g_{\operatorname*{can}}
\]
on $\left( -1,1\right) \times S^{n}$, which may be naturally identified with the
sphere $S^{n+1}$ with its north and south poles removed. The quantity $\psi(x,t)>0$
may thus be regarded as the radius of the hypersurface $\left\{ x\right\} \times
S^{n}$ at time $t$. It is natural to write geometric quantities related to $g$ in
terms of the distance $s\left( x\right) =\int_{0}^{x}\varphi(\hat{x})\,d\hat{x}$ from
the equator. Then one can write the metric in the nicer form of a warped product
\[
g=ds^{2}+\psi^{2}\,g_{\operatorname*{can}},
\]
bearing in mind that $s$ is ultimately a function of both $x$ and $t$.  (Note that
one encounters the commutator $[\partial_{t},\partial_{s}]=-n\psi
^{-1}\psi_{ss}\,\partial_{s}$ when taking derivatives with $x$ held constant.)

In \cite{AK}, we established neckpinching for a class of data essentially described
by three conditions: (1) the initial metric should have positive scalar curvature on
all of $S^{n+1}$ and positive Ricci curvature on the polar caps (in the terminology
of \cite{AK}); (2) its sectional curvature should be positive on planes tangential to
the spheres $\{x\}\times S^{n}$; and (3) it should be `sufficiently pinched',
i.e.~the minimum radius should be sufficiently small relative to the maximum
radius. In this paper, we impose an additional hypothesis: (4) the initial metric
should be reflection symmetric, i.e.~$\psi(s,0)\equiv\psi(-s,0)$. Under these
hypotheses, the results we obtain in Section \ref{Round} (combined with results from
the predecessor \cite{AK} to this paper) imply the the solution obeys a precise
asymptotic profile.

To describe this profile requires some additional notation. If $T<\infty$ is the
singularity time, let $u=\psi/\sqrt{2(n-1)(T-t)}$ be the blown-up radius defined below
in (\ref{eq:u-def}); let $\sigma=s/\sqrt{T-t}$ be the rescaled distance to the
neck defined in (\ref{eq:sigma-def}); and let $\tau=-\log(T-t)$ denote the
rescaled time variable defined in (\ref{eq:tau-def}). Note that $(\sigma,\tau)$
are self-similar coordinates with respect to the blown-up flow. We summarize our
results as

\begin{theorem}
For an open set of initial metrics symmetric with respect to reflection and
rotation, the solution $(S^{n+1},g(t))$ of Ricci flow becomes singular at
$T<\infty$. Its diameter remains bounded for all $t\in\lbrack0,T)$. The
singularity occurs only on the hypersurface $\left\{  0\right\}  \times S^{n}%
$. The solution satisfies the following asymptotic profile.

\textbf{Inner region:} on any interval $|\sigma|\leq A$, one has
\[
u(\sigma,\tau)=1+\frac{\sigma^{2}-2}{8\tau}+o(\frac{1}{\tau})\qquad
\text{uniformly as }\tau\rightarrow\infty.
\]

\textbf{Intermediate region:} on any interval $A\leq|\sigma|\leq B\sqrt{\tau}%
$, one has
\[
u(\sigma,\tau)=\sqrt{1+\left(  1+o(1)\right)  \frac{\sigma^{2}}{4\tau}}%
\qquad\text{uniformly as }\tau\rightarrow\infty.
\]

\textbf{Outer region:} for any $\varepsilon>0$, there exist $C<\infty$ and
$\bar{t}<T$ such that%
\[
\left(  \frac{1}{2}\sqrt{n-1}-\varepsilon\right)  \frac{s}{\sqrt{\log(1/s)}%
}\leq\psi(s,t)\leq\left(  \frac{1}{2}\sqrt{n-1}+\varepsilon\right)  \frac
{s}{\sqrt{\log(1/s)}}%
\]
for all $|\sigma|\geq C\sqrt{\tau}$ and $t\in(\bar{t},T)$.
\end{theorem}

\medskip

\noindent\textbf{Formal asymptotics. }In Section \ref{NotRound}, we derive
formal matched asymptotics for fully general neckpinches, without any symmetry
assumptions whatsoever. These are analogous to the formal asymptotics for the $m=2$
(presumably typical) case of \textsc{mcf} singularities considered in \cite{AV}. Our
method here is to study arbitrary perturbations of the self-similar cylinder
soliton. That is, we consider $\tilde{g}=g+h$, where
\[
g=dx^{2}+2(n-1)(T-t)\,g_{\operatorname*{can}}
\]
is the cylinder soliton on $\mathbb{R}\times S^{n}$, and $h$ is an arbitrary (small)
$(2,0)$-tensor. To accomplish our analysis, we modify the Ricci flow so that the
soliton becomes a fixed point of a related strictly parabolic flow. Because the
linearization of that flow has a null eigenvalue, one expects to find a center
manifold. We proceed by carrying out a quadratic variational analysis in order
formally to compute the dynamics of the flow on that center manifold.

The computations involved are extensive, but the conclusion we obtain is quite
satisfying: the formal matched asymptotics suggest that the behavior analyzed
rigorously in Section \ref{Round} should indeed be stable for fully general Ricci
flow neckpinches.

\subsection{A diameter bound}

We conclude this introduction with an observation that follows directly from the
estimates obtained in \cite{AK}. This result confirms an expectation of G.~Perelman
\cite{Hp}.

\begin{lemma}\label{TheNapkinLemma}
  Let $(S^{n+1},g(t))$ be any $\operatorname*{SO}(n+1)$-invariant solution of Ricci
  flow such that $g(0)$ has positive scalar curvature and positive sectional
  curvature on planes tangential to the spheres $\{x\}\times S^{n}$.

  Assume the metric $g(t)$ has at least two bumps for all $t<T$ (in the language of
  \cite{AK}).  Let $x=a(t)$ and $x=b(t)$ be the locations of the left- and right-most
  bumps, and assume that for all $t<T$ one has $\psi(a(t),t)\geq c$ , $\psi(b(t),
  t)\geq c$ for some constant $c>0$.

  If $g(t)$ becomes singular at $T<\infty$, then $\operatorname*{diam}
  (S^{n+1},g(t))$ remains bounded as $t\nearrow T$.
\end{lemma}

\begin{proof}
  By Proposition 5.4 of \cite{AK}, the limit profile $\psi(\cdot,T)$ exists.  Let
  $a(t)<b(t)$ in \thinspace$(-1,1)$ denote the left-most and right-most bumps,
  respectively. Then $a(t)\rightarrow a(T)$ and $b(t)\rightarrow b(T)$.  By Lemma
  5.6, the Ricci curvature is positive (so distances are decreasing) on $(-1,a(t)]$
  and $[b(t),1)$. Hence it will suffice to bound
  $\operatorname*{dist}{}_{g(t)}(x_{1},x_{2})$ for arbitrary $x_{1}<x_{2}$ in
  $(a(T)-\varepsilon,b(T)+\varepsilon)\subset(-1,1)$.

  Equations (5) and (11) of \cite{AK} imply that
  \begin{align*}
    \frac{d}{dt}\operatorname*{dist}{}_{g(t)}(x_{1},x_{2})  &  =\frac{d}{dt} 
    \int_{x_{1}}^{x_{2}}\varphi(x,t)\,dx\\
    &  =n\int_{s(x_{1})}^{s(x_{2})}\frac{\psi_{ss}}{\psi}\,ds\\
    &  =n\left\{ 
             \left. \frac{\psi_{s}}{\psi} 
             \right\vert _{s(x_{1})}^{s(x_{2})}
        +\int_{s(x_{1})}^{s(x_{2})}
           \left(  \frac{\psi_{s}}{\psi}\right)^{2}\,ds
           \right\}  .
  \end{align*}
  Proposition 5.1 of \cite{AK} bounds $\psi_{s}$ uniformly, while Lemma 5.5 shows
  that the number of bumps and necks is nonincreasing in time. It follows that
  \[
  \int_{s(x_{1})}^{s(x_{2})}\left( \frac{\psi_{s}}{\psi}\right)
  ^{2}\,ds\leq
  C\int_{s(x_{1})}^{s(x_{2})}\frac{|\psi_{s}|}{\psi^{2}}\,ds\leq
  C[\frac{1}
  {\psi_{\min}(t)}-\frac{1}{\psi_{\max}(t)}]\leq\frac{C}{\psi_{\min}(t)}.
  \]
  Hence Lemma 6.1 lets us conclude that
  \[
  \left\vert
    \frac{d}{dt}\operatorname*{dist}{}_{g(t)}(x_{1},x_{2})\right\vert
  \leq\frac{C}{\psi_{\min}(t)}\leq\frac{C}{\sqrt{T-t}},
  \]
  which is integrable.
\end{proof}

Recalling Lemmata 7.1 and 10.1 of \cite{AK}, one immediately obtains

\begin{proposition}
  Assume that $(S^{n+1},g(t))$ satisfies the hypotheses of Lemma \ref{TheNapkinLemma}
  and is reflection-symmetric. Then the final-time profile $\psi(\cdot,T)$ has
  $\psi(x,T)>0$ for all $x\neq0$. Because $|\operatorname*{Rm}|\leq C/\psi^{2}$, this
  implies that the singularity occurs only on the hypersurface $\{0\}\times S^{n}$.
\end{proposition}

\medskip

Diameter estimates for arbitrary finite-time Ricci flow singularities are studied
from another point of view in \cite{T}. Lower bounds for diameter are known for
certain topologies; see \cite{IK}.

\section{Rigorous asymptotics of the $\mathrm{SO}(\R^{n+1})$ invariant neckpinch}
\label{Round}

\subsection{The blown-up radius and the linearized equation}
The radius $\psi$ satisfies
\begin{equation}
  \label{eq:psi-evolution}
  \psi_t = \psi_{ss}-(n-1)\frac{1-\psi_s^2}{\psi}.
\end{equation}
We will assume that the solution is defined on some time interval $t_0\leq t<T$,
where the initial time $t_0$ neeed not be 0. In order to construct examples of
solutions that follow our precise asymptotics, we will at some point actually have to
assume that $T-t_0$ is ``sufficiently small''.

We consider the blown-up radius
\begin{equation}
  \label{eq:u-def}
  u \doteqdot \frac\psi{R_n(t)},
  \quad \text { where } \quad
  R_n(t)\doteqdot \sqrt{2(n-1)(T-t)},
\end{equation}
in which $T$ may or may not be the blow-up time. We also introduce the rescaled
distance to the neck
\begin{equation}
  \label{eq:sigma-def}
  \sigma\doteqdot \frac s {\sqrt{T-t}}.
\end{equation}
(Here the absence of a factor $n-1$ in $\sqrt{T-t}$ is intentional and will result in
the numerically simplest equations later on.) Then we have
\[
\psi_s= \sqrt{2 (n-1)}\;u_\sigma,\qquad \psi_{ss}=\sqrt{2
  (n-1)}\,\frac{u_{\sigma\sigma}}{\sqrt{T-t}},
\]
and
\[
\psi_t = \sqrt{2 (n-1) (T-t)}\; u_t - \frac{\sqrt{2 (n-1)}}
{2\sqrt{T-t}}\;u.
\]
This leads to the following evolution equation for $u$,
\begin{equation}
  \label{eq:u-evolution}
  u_\tau = u_{\sigma\sigma}
  +\frac12\l(u-\frac1u\r)
  +(n-1)\frac{u_\sigma^2}{u},
\end{equation}
where
\begin{equation}
  \label{eq:tau-def}
  \tau \doteqdot -\log(T-t).
\end{equation}
The blown-up time variable $\tau$ takes values in the interval
\[
\tau_0\leq \tau<\infty, \text{ where } \tau_0 \doteqdot -\log (T-t_0).
\]
Without loss of generality we will assume that $\tau_0\geq 1$.

Near the neck, $u$ will be close to the self-similar radius $u\approx1$. So to
linearize, we set
\begin{equation}
  \label{eq:v-def}
  u=1+v
\end{equation}
and observe that $v$ must satisfy
\begin{equation}
  \label{eq:v-evolution}
  v_\tau= v_{\sigma\sigma}+v+
  \frac{(n-1)v_\sigma^2-\frac12v^2}
  {1+v}.
\end{equation}
Here the partial derivatives $\pd/\pd\sigma$ and $\pd/\pd\tau$ do not commute. One
has instead
\[
[\pd_t, \pd_s] = -n \frac{\psi_{ss}}{\psi} \pd_s,
\]
and thus
\begin{equation}
  \label{eq:commutator}
  [\pd_\tau,\pd_\sigma] =
  \l[(T-t)\pd_t , \sqrt{T-t}\,\pd_s\r] =
  -\l\{\frac12 + n\frac{u_{\sigma\sigma}}{u}\r\}\pd_\sigma.
\end{equation}

\subsection{Pointwise estimates for $u$ and $v$}
In \cite{AK}, we showed that on solutions of the Ricci flow whose
scalar curvature is nonnegative, the quantity
\[
F = \frac KL\l\{2+\log \frac L {L_{\mathrm{min}}(0)}\r\}
\]
satisfies a maximum principle. Here
\[
K=\frac{\psi_{ss}}{\psi} \text{ and } L=\frac{1-\psi_s^2}{\psi^2},
\]
and
\[
L_{\mathrm{min}}(0) = \inf \l. \frac{1-\psi_s^2}{\psi^2}\r|_{t=0}.
\]
The maximum principle implies that $\sup F$ does not increase when it
is above $n-1$.

Let us assume throughout this paper that the initial metric has
nonnegative scalar curvature, and that
\[
\sup F \leq F^*
\]
for some fixed constant $F^*\geqslant n-1$.

When one writes the quantity $F$ in terms of $\psi$, $\psi_s$, and
$\psi_{ss}$, the estimate $F\leq F^*$ implies the following upper
bounds for $\psi$, $\psi_s$, and $\psi_{ss}$ in terms of the radius of
the neck $\psi(0,t)$.
\begin{lemma}\label{lem:pointwise-u-estimates}
  There exist constants $\delta>0$ and $C<\infty$ that only depend on
  $F^*$ and the dimension $n$ such that
  \begin{equation}
    \label{eq:u-inner-estimate}
    1-\frac C\tau
    \leq u
    \leq 1+C\frac{1+\sigma^2}{\tau}
    \text{ for }|\sigma|\leq2\surd \tau,
  \end{equation}
  and
  \begin{equation}
    \label{eq:u-intermediate-estimate}
    1-\frac C\tau\leq u\leq
    C\frac {|\sigma|}{\surd{\tau}} \sqrt{\log\frac {|\sigma|}{\surd{\tau}}}
    \text{ for }
    2\surd\tau\leq |\sigma|\leq e^{\delta\tau}.
  \end{equation}
  For the derivative $u_\sigma$, we have
  \begin{gather}
    |u_\sigma| \leq C\frac{1+|\sigma|}{\tau}
    \text{ for }|\sigma|\leq2\surd \tau, \label{eq:usigma-inner} \\
    |u_\sigma| \leq \frac C{\surd\tau}\sqrt{{\log\frac
        {|\sigma|}{\surd{\tau}}}} \text{ for } 2\surd\tau\leq
    |\sigma|\leq e^{\delta\tau}.\label{eq:usigma-intermediate}
  \end{gather}
\end{lemma}
In both cases, these estimates imply that $|u_\sigma|\leq C\delta$ at
all points with $|\sigma|\leq e^{\delta\tau}$. After decreasing
$\delta$, if necessary, we may assume that $C\delta$ is as small as we
like. In particular, in the region we are considering, we always may
assume that
\begin{equation}
  \label{eq:usigma-small}
  |u_\sigma|\leq \frac12.
\end{equation}
Since $v=u-1$, we immediately get estimates for $v$ and $v_\sigma$.
Those for $v_\sigma$ are of course the same as for $u_\sigma$. For
$v$, we have
\begin{equation}
  \label{eq:v-pointwise}
  \l\{
  \begin{split}
    &-\frac C{\tau} \le v \le C\frac{1+\sigma^2}{\tau}
    \text{ for }|\sigma|\le 2\surd\tau,\\
    &-\frac C{\tau} \le v \le C \frac{|\sigma|}{\surd\tau}
    \sqrt{\log\frac{|\sigma|}{\surd\tau}} \text{ for }2\surd\tau\le
    |\sigma|\le e^{\delta\tau},
  \end{split}
  \r.
\end{equation}
in which $C$ is again a constant depending only on $F^*$ and $n$.

\begin{proof}[Proof of Lemma \ref{lem:pointwise-u-estimates}]
  The upper estimates follow immediately from Lemma 9.4 in \cite{AK}.
  Only the lower estimate in (\ref{eq:u-inner-estimate}) was not
  proved in \cite{AK}. Briefly, the lower estimate follows from the
  boundedness of $F$. At the neck, one has
  \[
  K=\frac{\psi_{ss}}\psi,\quad L=\frac1{\psi^2},\quad
  \psi=(1+o(1))\sqrt{2(n-1)(T-t)},
  \]
  so that $(K/L)\log L\leq C$ implies
  \[
  0\leq \psi_{ss}\leq \frac{C}{\tau\psi}.
  \]
  The radius $\psi(0,t)$ of the neck satisfies
  \[
  \psi_t = \psi_{ss}-\frac{n-1}{\psi} \leq -\frac{n-1 - C
    (\tau^{-1})}{\psi},
  \]
  which upon integration yields
  \[
  \psi(0,t) \geq \l(1-\frac{\tilde C}\tau\r)\sqrt{2 (n-1) (T-t)}
  \]
  where $\tilde C$ only depends on $F^*$ and $n$.  This implies the
  lower estimates in the Lemma.
\end{proof}

\subsection{Transition to commuting variables}
\label{sec:commuting-vars}
On most of the neck, the quantity $v=o(1)$. So it is natural to drop
the quadratic terms in \eqref{eq:v-evolution}, resulting in the
deceptively simple equation
\[
v_\tau= v_{\sigma\sigma}+v
\]
for $v$. Here the partial derivatives $\pd_\tau$ and $\pd_\sigma$ do
not commute. To work with commuting variables, we regard $v$ as
function of $\sigma$. In other words, we define
\begin{equation}
  \label{eq:sigma}
  \sigma(x,t) = \frac1{\sqrt{T-t}}\int_0^x ds 
  = e^{-\tau/2}\int_0^x ds, 
\end{equation}
and consider a function $\tilde v(\sigma,\tau)$ such that $\tilde
v(\sigma(x,\tau),\tau)\equiv v(x,\tau)$. We shall abuse notation and
write $v$ for both quantities $v$ and $\tilde v$. This is only
ambiguous when we take time derivatives, in which case we write
\[
\l.  \frac{\pd v}{\pd \tau}\r|_{\sigma\text{ const}} = \frac{\pd
  \tilde v}{\pd \tau} \quad\text{and}\quad \l.  \frac{\pd v}{\pd
  \tau}\r|_{x\text{ const}}= \frac{\pd v}{\pd \tau}.
\]
The two are related by
\[
\l.\frac{\pd v}{\pd \tau}\r|_{\sigma\text{ const}} =\l.\frac{\pd
  v}{\pd \tau}\r|_{x \text{ const}} -\frac{\pd \sigma}{\pd \tau}
\frac{\pd v}{\pd\sigma}.
\]
\begin{lemma}\label{lem:J}
  Assuming reflection symmetry of the metric, i.e.\ $\psi(-s,t)\equiv
  \psi(s,t)$, one has
  \[
  \frac{\pd \sigma}{\pd \tau} = \frac\sigma2+nJ(\sigma,\tau),
  \]
  where
  \[
  J(\sigma,\tau) = \int _0^\sigma
  \frac{v_{\sigma\sigma}}{1+v}\;d\sigma = \frac{v_\sigma}{1+v}
  +\int_0^\sigma \frac{v_\sigma^2}{(1+v)^2}\;d\sigma.
  \]
\end{lemma}
\begin{proof}
  By definition one has $\sigma(x, \tau) = e^{\tau/2}s(x,
  T-e^{-\tau})$. Hence
  \[
  \sigma_\tau = \frac12\sigma + e^{\tau/2}s_t(x, T-e^{-\tau})e^{-\tau}
  = \frac12\sigma + e^{-\tau/2}s_t(x, T-e^{-\tau}).
  \]
  One also has \cite[\S10]{AK}
  \[
  \frac{\partial s}{\partial t} =\int_0^s \frac{n\psi_{ss}}{\psi} \,
  ds =n\l\{\frac{\psi_s}{\psi} + \int_0^s
  \frac{\psi_s^2}{\psi^2}\,ds\r\}.
  \]
  Using the relations (\ref{eq:u-def}) and (\ref{eq:sigma-def})
  between $u=1+v$ and $\psi$, and $s$ and $\sigma$, one finds the
  stated expression for $\partial_\tau\sigma$.
\end{proof}

From here on we will consider $v$ as a function of $\sigma$ and
$\tau$. All $\tau$ derivatives are intended to be time derivatives
with $\sigma$ kept constant. It follows from \eqref{eq:v-evolution}
that $v$ satisfies
\begin{equation}
  \label{eq:v-sigma}
  v_\tau = v_{\sigma\sigma}-\frac\sigma2 v_\sigma + v
  -nJ(\sigma,\tau) v_\sigma
  +\frac{(n-1)v_\sigma^2-\frac12v^2}
  {1+v}.
\end{equation}
We write this equation as
\[
v_\tau= A(v) + N(v),
\]
where $A$ is the linear differential operator
\[
A = \l(\frac\pd{\pd\sigma}\r)^2 - \frac\sigma2 \frac\pd{\pd\sigma}+1,
\]
and $N(v)$ represents the remaining (nonlinear) terms in
\eqref{eq:v-sigma},
\[
N(v) = -nv_\sigma\int_0^\sigma\frac{v_{\sigma\sigma}}{1+v}d\sigma
+\frac{(n-1)v_\sigma^2-\frac12v^2}{1+v}.
\]
Integrate by parts to get
\begin{equation}\label{eq:Nv-form1}
  N(v) = -nv_\sigma\int_0^\sigma\frac{v_\sigma^2}{(1+v)^2}\,d\sigma
  -\frac{v_\sigma^2+\frac12v^2}{1+v}.
\end{equation}

\begin{lemma}\label{lem:rough-Nv-estimate}
  For all $|\sigma|\leq e^{\delta\tau}$, one has
  \[
  |N(v)|\leq C\frac{1+\sigma^4}{\tau^2}.
  \]
\end{lemma}
\begin{proof}
  This follows from Lemma \ref{lem:pointwise-u-estimates}. In fact,
  Lemma~\ref{lem:pointwise-u-estimates} implies that
  \begin{equation}\label{eq:1600}
    |v|\leq C\frac{1+\sigma^2}{\tau}
    \text{ and }
    |v_\sigma|\leq C\frac{1+|\sigma|}{\tau}
  \end{equation}
  for $|\sigma|\leq 2\surd\tau$. For larger $\sigma$, one finds other
  estimates for $v$ and $v_\sigma$, which are stronger than
  \eqref{eq:1600}. Using \eqref{eq:1600} for all $|\sigma|\leq
  e^{\delta\tau}$, one arrives at our estimate for $N(v)$.
\end{proof}

\subsection{The linearized equation} \label{sec:linearized}
If we ignore the nonlinear terms in \eqref{eq:v-sigma}, then we see
that the small quantity $v$ satisfies $v_\tau=A(v)$. Even though $v$
is not defined for \emph{all} $s\in\R$, we know from our estimate on
the neck that $v$ is defined for $|\sigma|\leq e^{\delta\tau}$, for
some $\delta=\delta(F^*,n)$. In the first analysis, we assume that $v$
is in fact a solution of the Cauchy problem $v_\tau=A(v)$. The
operator $A$ is self-adjoint in the Hilbert space
\[
\fh=\l\{v\in L^2\l(\R, e^{-\sigma^2/4}d\sigma\r) \mid v(\sigma)\equiv
v(-\sigma) \r\}.
\]
It has pure point spectrum, with eigenvalues $\lambda_m=1-m$ for
$m=0,1,2,\ldots$. The eigenfunction corresponding to the eigenvalue
$\lambda_m$ is an even \emph{Hermite polynomial} of degree $2m$.
Normalizing the eigenfunction so that the coefficient of its highest
order term is 1, we set
\[
h_{2m}(\sigma) = \sigma^{2m} -\frac{(2m)!}{1!(2m-2)!}\sigma^{2m-2}
+\frac{(2m)!}{2!(2m-4)!}\sigma^{2m-4} -\cdots
\]
In particular,
\begin{equation}
  \label{eq:h024-explicit}
  h_0(\sigma)=1, \quad
  h_2(\sigma)=\sigma^2-2,\quad
  h_4(\sigma)=\sigma^4-12\sigma^2+12.
\end{equation}

\subsection{Eigenfunction decomposition of $v$}
Let $\eta\in C^\infty(\R)$ be an even bump function with $\eta(z)=0$
for $z\geq 1\frac14$ and $\eta(z)=1$ for $z\leq1$. We then define
\begin{equation}
  \label{eq:V-def}
  V(\sigma,\tau) = \eta(\sigma e^{-\delta\tau/2})v(\sigma,\tau).
\end{equation}
For $|\sigma| \geq \frac54 e^{\delta\tau/2}$, we set $V\equiv0$. Then
$V(\cdot,\tau)\in\fh$, and $v\equiv V$ for $|\sigma|\leq
e^{\delta\tau/2}$.

A computation shows that
\begin{equation}
  \label{eq:V-evolution}
  V_\tau-A(V) = \eta N(v)+E,
\end{equation}
where the ``error term'' is
\[
E=\l(\eta_\tau-\eta_{\sigma\sigma}+\frac\sigma2\eta_\sigma\r)v
-2\eta_\sigma v_\sigma.
\]
\begin{lemma}\label{lem:E-estimate}
  The error term $E$ vanishes except when $e^{\delta\tau/2} \leq
  |\sigma|\leq \frac54 e^{\delta\tau/2}$.  When $E\neq 0$, one has
  \[
  |E(\sigma,\tau)|\leq C|\sigma|,
  \]
  where the constant $C$ only depends on $F^*$ and $n$.  One also has
  \[
  \l\| E (\cdot,\tau)\r\| \leq C \exp\l(-e^{\delta\tau/2}\r),
  \]
  where the constant $C$ again only depends on $F^*$ and $n$.
\end{lemma}
\begin{proof}
  The pointwise estimate follows from the boundedness of $\eta_\tau$,
  $\sigma\eta_\sigma$, and $\eta_{\sigma\sigma}$, as well as the
  pointwise bounds for $v$ provided by Lemma
  \ref{lem:pointwise-u-estimates}.  Given the pointwise bounds for
  $E$, one finds that
  \[
  \|E\| ^2 \leq C^2\int_{e^{\delta\tau/2}}^\infty
  \sigma^2e^{-\sigma^2/4}\,d\sigma,
  \]
  from which the $\fh$ norm estimate follows.
\end{proof}

For each $\tau$, the function $V(\cdot,\tau)$ belongs to $\fh$. So we
can consider the splitting into mutually orthogonal terms determined
by
\begin{equation}
  \label{eq:V-decomposition}
  V(\sigma,\tau) = a_0(\tau)h_0(\sigma) + a_2(\tau)h_2 (\sigma) + W(\sigma,\tau).
\end{equation}
Our pointwise estimate~\eqref{eq:v-pointwise} for $v$ suggests that
for large $\tau$ these terms will decay like $\tau^{-1}$. This turns
out to be the case for the middle term $a_2(\tau)h_2 (\tau)$, but the
other two terms are in fact smaller. We will show in the next few
paragraphs that they decay like $\tau^{-2}$. This establishes
$a_2(\tau)h_2 (\tau)$ as the dominant term for large $\tau$.

\subsection{Easy estimates for $a_0$ and $a_2$}
Lemma 1 and in particular (\ref{eq:v-pointwise}) imply that
\[
|V(\sigma, \tau)|\leq \frac C\tau (1+\sigma^2)
\]
for all $\sigma\in\R$, and all $\tau\geq \tau_0$. Upon taking the
$\fh$-inner product with $h_0$ and $h_2$ this gives us
\begin{equation}\label{eq:a0-a2-bound}
  |a_0(\tau)|+|a_2(\tau)| \leq \frac C\tau
\end{equation}
for all $\tau\geq\tau_0$.

\subsection{Decay of $a_0$}
Since $\|h_0\| ^2 a_0(\tau)= (h_0, V)_\fh$, we have
\begin{align*}
  \|h_0\| ^2 a_0' (\tau)
  &= (h_{0}, V_{\tau})_{\fh} \\
  &= (h_{0}, A(V)+\eta N(v)+E)_{\fh}\\
  &= (h_0, V)_\fh+ (h_0, \eta N(v)+E)_\fh \\
  &= \|h_0\| ^2 a_0 (\tau) + (h_0, \eta N(v)+E)_\fh.
\end{align*}
Hence
\[
a_0'(\tau)-a_0 (\tau) = \l (\frac{h_0}{\|h_0\|^2}, \eta N(v)+E\r)_\fh.
\]
Using our estimates for $N(v)$ and $E$, we conclude that
\[
f_0(\tau)\doteqdot a_0'(\tau)-a_0 (\tau)
\]
satisfies
\[
|f_0 (\tau) |\leq \frac{C(F^*,n)}{\tau^2}.
\]
The variation of constants formula tells us that
\[
a_0 (\tau) = -\int_\tau^\infty e^{\tau-\tau'}f_0 (\tau')\,d\tau',
\]
and thus that
\[
|a_0 (\tau)|\leq \frac{C(F^*,n)}{\tau^2}.
\]
\subsection{Decay of $W$}
Since $W$ lies in the stable space of the operator $A$, i.e.\ the
space spanned by those eigenfunctions with negative eigenvalues, we
expect $W$ to decay according to the slowest stable eigenvalue
(i.e.~$\sim e^{-\tau}$), or else to decay like the ``forcing term''
$\eta N(v) +E$. In this situation, the forcing term is dominant, and
we have the following.
\begin{lemma}\label{lem:W-decay}
  For any given $\tau_0>0$, one has for all $\tau>\tau_0$
  \[
  \|W(\cdot,\tau)\| \le e^{\tau_0-\tau} \|W(\cdot,\tau_0)\| + \frac
  C{\tau^2}.
  \]
  Consequently, one also has
  \[
  \|W(\cdot, \tau)\| \leq \frac C{\tau^2} \bigl (1+W_0 \bigr).
  \]
  Here the constant $C$ only depends on $F^*$ and $n$, and
  \[
  W_0 \doteqdot e^{\tau_0}\|W(\cdot,\tau_0)\| .
  \]
\end{lemma}

\begin{proof}
  Since $\|W\|^2=(W,W)=(W,V)$ and $(W,A(V))_\fh=(W,A(W))_\fh$, one has
  \begin{align*}
    \|W\| \;\frac{d\|W\| }{d\tau}
    &=\frac12 \frac d{d\tau}\|W\| ^2 \\
    &= \frac12 \frac d{d\tau} \l (W, V\r )_\fh \\
    &= \l (W, V_\tau\r)_\fh \\
    &= \l (W, A(V)+\eta N(v)+E\r)_\fh \\
    &= (W,A(W))_\fh + \l(W, \eta N(v)+E\r)_\fh \\
    &\leq -\|W\| ^2 +\|W\|\cdot\| \eta N(v)+E\| ,
  \end{align*}
  whence
  \begin{equation}
    \label{eq:W-decay1}
    \frac{d\|W\| }{d\tau} \le -\|W\| +\| \eta N(v)+E\| .
  \end{equation}
  From Lemma~\ref{lem:rough-Nv-estimate} we conclude that $|\eta
  N(v)|\le C (1+\sigma^4)\tau^{-2}$ for all $\sigma$, and hence that
  \[
  \| \eta N(v)\| \le \frac C{\tau^2}.
  \]
  For $E$, we use Lemma~\ref{lem:E-estimate} together with the
  calculus inequality $Ce^{-e^{\delta\tau/2}}\le C'\tau^{-2}$. We
  apply all this to \eqref{eq:W-decay1} and conclude that
  \begin{equation*}
    \frac{d\|W\| }{d\tau} \le -\|W\|  + \frac C{\tau^2}.
  \end{equation*}
  Therefore
  \[
  \frac{d(e^\tau\|W\| )}{d\tau} \le Ce^\tau \tau^{-2}.
  \]
  Integration leads to
  \begin{equation*}
    e^\tau \|W(\cdot,\tau)\| 
    \le e^{\tau_0} \|W(\cdot, \tau_0)\|  
    +C \int_{\tau_0}^\tau \frac{e^{\tilde\tau}}{\tilde\tau^2}\;d\tilde\tau
    \le e^{\tau_0} \|W(\cdot, \tau_0)\|  + \frac{\tilde C}{\tau^2}e^\tau,
  \end{equation*}
  as claimed. The second estimate for $\|W\| $ follows directly from
  $\sup \tau^2 e^{-\tau}<\infty$.
\end{proof}

\subsection{Derivative estimates for $W$}
We use the regularizing effect of the heat equation to bootstrap the
estimates of Lemma~\ref{lem:W-decay} by one space derivative.
\begin{lemma}\label{eq:W-deriv-est}
  Let $\tau_0>0$ be as in Lemma~\ref{lem:W-decay}. Then for all
  $\tau>\tau_0+1$, one has
  \[
  \|W_\sigma(\cdot,\tau)\| \le C(1+W_0)\tau^{-2},
  \]
  where the constant $C$ only depends on $F^*$ and $n$, and as before,
  $W_0 = e^{\tau_0}\|W(\cdot,\tau_0)\| $.
\end{lemma}
\begin{proof}
  Note that
  \[
  Au = u_{\sigma\sigma}-\frac12\sigma u_\sigma + u =
  e^{\sigma^2/4}\bigl( e^{-\sigma^2/4}u_\sigma\bigr)_\sigma
  \]
  so that for $u\in D(A)$ one has
  \[
  -(u,Au)_\fh = \int_\R e^{-\sigma^2/4}
  \bigl\{u_\sigma^2-u^2\bigr\}\,d\sigma.
  \]
  Thus
  \begin{equation}
    \label{eq:u-deriv-bound}
    \|\pd_\sigma u\| ^2 = \|u\| ^2 -(u,Au).
  \end{equation}
  If $u\perp h_0,h_2$ then $(u,-Au)_\fh\geq0$, and we have
  \[
  \|\pd_\sigma u\| \leq \|u\| +\surd (u, -Au)_\fh.
  \]
  Applying the spectral theorem to the selfadjoint operator $A$ on the
  Hilbert space $\{h_0,h_2\}^\perp$ in order to estimate $(u,
  Ae^{\theta A}u)$, one concludes that
  \begin{equation}
    \label{eq:semigroup-deriv-est}
    \left\|\pd_\sigma e^{\theta A}u\right\|  \leq\frac C{\surd\theta}\|u\| 
  \end{equation}
  holds for $0<\theta\leq 1$, and all $u\in D(A)$.

  The stable component $W$ of $V$ satisfies
  \[
  W_\tau = A(W) +f(\tau), \text{ where }f(\tau)\doteqdot\cP(\eta
  N(v)+E),
  \]
  and where $\cP$ is the $\fh$-orthogonal projection onto the space
  $\{u\in\fh\mid u\perp h_0,h_2\}$.  By the variation of constants
  formula one then has, for all $\tau\geq\tau_0+1$
  \[
  W(\tau) = e^AW(\tau-1) + \int_0^1 e^{\theta
    A}f(\tau-\theta)\;d\theta.
  \]
  Since $f(\tau)\perp \{h_0,h_2\}$, we can use
  (\ref{eq:semigroup-deriv-est}) to get
  \begin{align*}
    \|\pd_\sigma W(\tau)\| &\leq \|\pd_\sigma e^AW(\tau-1)\| +
    \int_0^1 \left\|\pd_\sigma e^{\theta A}f(\tau-\theta)\right\|
    \;d\theta \\
    &\leq C\|W(\tau-1)\| + \int_0^1 \frac C{\surd\theta}
    \left\|f(\tau-\theta)\right\| \;d\theta.
  \end{align*}
  The combination of Lemmas \ref{lem:rough-Nv-estimate} and
  \ref{lem:E-estimate} again gives $\|f (\tau)\| \leq C\tau^{-2}$.
  The integral can therefore be bounded by $C\tau^{-2}$. Combining
  this with our previous estimate for $\|W (\tau-1)\| $ from Lemma
  \ref{lem:W-decay}, we get the inequality in Lemma
  \ref{eq:W-deriv-est}.
\end{proof}

\subsection{A bound for $\|\sigma W\| $}
If $\varphi\in C^1_c(\R)$, then one has the identity
\[
\int_\R e^{-\sigma^2/4}2\varphi(\sigma)\sigma\varphi'(\sigma)\,d\sigma
=\int_\R e^{-\sigma^2/4} \bigl( -2\varphi(\sigma)^2 +
\tfrac12\sigma^2\varphi(\sigma)^2 \bigr) \,d\sigma
\]
which after rearranging leads to
\begin{align*}
  \|\sigma \varphi\| ^2 &=4\int_R e^{-\sigma^2/4}
  (\sigma\varphi\varphi'+\varphi^2)\,d\sigma \\
  &=4\l\{(\sigma\varphi, \varphi')_\fh + \|\varphi\| ^2 \r\} \\
  &\le \tfrac12 \|\sigma\varphi\| ^2 + 8 \|\varphi'\| ^2
  + 4\|\varphi\| ^2 \\
  &\le 16\l\{\|\varphi'\| ^2 + \|\varphi\| ^2\r\}.
\end{align*}
An approximation argument leads to
\begin{lemma}\label{lem:sigmaW-bound}
  For any $\varphi\in\fh$ with $\varphi_\sigma\in\fh$ one has
  $\sigma\varphi\in\fh$, while
  \[
  \|\sigma\varphi(\sigma)\| \leq 4\l\{\|\varphi_\sigma\| + \|\varphi\|
  \r\}.
  \]
  In particular, $W$ satisfies
  \begin{equation}
    \label{eq:sigmaW-estimate}
    \|\sigma W\|  
    \leq  C(1+W_0)\tau^{-2}
  \end{equation}
  for all $\tau\ge\tau_0+1$.
\end{lemma}

\subsection{Pointwise estimates for $W$}
A Sobolev inequality says that the bound on $\|\pd_\sigma W\| $
implies a pointwise bound for $W$. Due to the exponential weight in
the norm this bound is not uniform. One has
\begin{lemma}\label{lem:Gauss-Sobolev}
  Any $\varphi\in\fh$ with $\pd_\sigma\varphi\in\fh$ is a continuous
  function, which satisfies
  \[
  e^{-\sigma^2/8} \l| \varphi(\sigma)\r| \leq C\bigl (\|\varphi\|
  +\|\pd_\sigma\varphi\| \bigr).
  \]
  Hence $W$ satisfies
  \[
  |W(\sigma,\tau)| \le \frac {C(1+W_0)}{\tau^2}e^{\sigma^2/8}.
  \]
\end{lemma}
\begin{proof}
  Assuming $\varphi\in C^1_c(\R)$ one has
  \begin{align*}
    e^{-\sigma^2/4}\varphi(\sigma)^2 & = \int_{-\infty}^\sigma
    e^{-\sigma^2/4} \bigl[ -\frac\sigma2\varphi^2 +
    2\varphi\varphi_\sigma\bigr]
    \,d\sigma \\
    &\le \|\sigma\varphi\| \|\varphi\|
    + \|\varphi\|  \|\varphi_\sigma\|  \\
    &\le C\bigr\{ \|\varphi\| ^2 + \|\varphi_\sigma\| ^2\bigl\}.
  \end{align*}
  So the pointwise estimate holds for $\varphi\in C^1_c$. By
  approximation it also holds for all $\varphi\in\fh$ with
  $\pd_\sigma\varphi\in\fh$.
\end{proof}

Lemma \ref{lem:Gauss-Sobolev} says that $W$ decays pointwise like
$\tau^{-2}$, and uniformly on bounded $\sigma$ intervals. The
following pointwise bound for $W$ only gives a $\tau^{-1}$ decay rate,
but it is stronger for large values of $\sigma$. The estimate also
does not depend on $W_0$.
\begin{lemma}\label{lem:W-pointwise-2}
  For some constant $C=C(F^*,n)<\infty$ one has for all
  $\tau>\tau_0+1$
  \[
  |W(\sigma,\tau)| \leq \frac C\tau (1+\sigma^2).
  \]
\end{lemma}
\begin{proof}
  We have, by definition, $W(\sigma,\tau) =
  V(\sigma,\tau)-a_0(\tau)h_0(\sigma)-a_2(\tau)h_2(\sigma)$. Using our
  pointwise bounds for $V$, $a_0$, $a_2$, and $h_0(\sigma)\equiv 1$,
  $h_2(\sigma)\equiv \sigma^2-2$, one quickly gets the stated estimate
  for $W$.
\end{proof}

\subsection{Pointwise estimate for $W_\sigma$}
By definition we have $W(\tau,\sigma) =
V(\tau,\sigma)-a_0(\tau)h_0(\sigma) - a_2(\tau)h_2(\sigma)$.
Differentiate and keep in mind that $h_0'(\sigma)=0$ and
$h_2'(\sigma)=2\sigma$. We get $ W_\sigma = V_\sigma -2a_2(\tau)\sigma
$. Using the pointwise derivative bounds (\ref{eq:usigma-inner}) and
(\ref{eq:usigma-intermediate}) for $u$, one then gets
\begin{equation}
  \label{eq:Wsigma-bounds}
  \l| W_\sigma\r| \leq \frac C\tau \bigl(1+|\sigma|\bigr)
\end{equation}
which holds for all $\tau\geq\tau_0$ and $\sigma\in\R$.

\subsection{An equation for $a_2$}
As with $a_0$, we have
\begin{align*}
  \|h_2\| ^2 a_2' (\tau)
  &= (h_{2}, V_{\tau})_{\fh} \\
  &= (h_{2}, A(V)+\eta N(v)+E)_{\fh}\\
  &= (A(h_2), V)_\fh+ (h_2, \eta N(v)+E)_\fh \\
  &= (h_2, \eta N(v)+E)_\fh.
\end{align*}
We use~\eqref{eq:Nv-form1} and the identity $1/ (1+v) = 1-v/ (1+v)$ to
rewrite $N(v)$ as
\begin{equation}
  \label{eq:Nv-form2}
  N(v)= N_2(v)+N_3(v),\text{ where }\l\{
  \begin{gathered}
    N_2(v) = -\l(v_\sigma^2+\tfrac12v^2\r) \\
    N_3(v) =
    -nv_\sigma\int_0^\sigma\frac{v_\sigma^2}{(1+v)^2}\,d\sigma
    +\frac{vv_\sigma^2+\frac12v^3}{1+v}.
  \end{gathered}\r.
\end{equation}
Here, $N_2(v)$ is the purely quadratic part, while $N_3(v)$ contains
the cubic and higher order terms.
\begin{lemma}\label{lem:N3-estimate}
  There is a constant $C=C(F^*,n)$ such that for all $\tau\geq\tau_0$
  \[
  \|\eta N_3(v)\| \le \frac{C}{\tau^3}.
  \]
\end{lemma}
\begin{proof}
  The pointwise estimates for $v$ and $v_\sigma$ from
  \eqref{eq:v-pointwise} imply that
  \[
  |v|\leq C\frac{1+\sigma^2}{\tau},\quad |v_\sigma|\le
  C\frac{1+|\sigma|}{\tau}
  \]
  whenever $|\sigma|\le e^{\delta\tau}$. Hence we have
  \[
  |N_3(v)| \le C \frac {1+\sigma^6}{\tau^3} \text{ for } |\sigma|\le
  e^{\delta\tau}.
  \]
  (The term with $v^3$ contributes the highest power in $\sigma$.)

  This implies the Lemma.
\end{proof}

We continue with our computation of $a_2'(\tau)$. We have
\begin{align*}
  \|h_2\| ^2 a_2' (\tau)
  &= (h_2, \eta N(v)+E)_\fh \\
  &= (h_2, \eta N_2(v)+\eta N_3(v)+E)_\fh\\
  &= (h_2, \eta N_2(v))_\fh + (h_2, \eta N_3(v)+E)_\fh.
\end{align*}
The last term satisfies
\[
\l| (h_2, \eta N_3(v)+E)_\fh \r| \le \|h_2\| \l (\|\eta N_3(v)\|
+\|E\| \r) \le \frac C{\tau^3}.
\]
We would now like to replace the $v$ in the quadratic expression
$N_2(v)$ by $V$, and then by the dominant term
in~\eqref{eq:V-decomposition}, i.e.\ $a_2(\tau)h_2(\sigma)$. Before we
do this, we estimate the errors produced by these replacements.
\begin{lemma}\label{lem:N2-N2-estimate}
  There is a constant $C=C(F^*,n)$ such that for all $\tau\ge\tau_0$,
  one has
  \[
  \l\|\eta \l (N_2(v)-N_2(V)\r)\r\| \le C e^{-e^{\delta\tau/2}}.
  \]
\end{lemma}
In other words $\eta\l (N_2(v)-N_2(V)\r)$ satisfies the same
super-exponentially small estimate as the error term $E$. (See
Lemma~\ref{lem:E-estimate}.)
\begin{proof}
  Since $V=\eta\cdot v$, a direct calculation gives
  \[
  \eta\l (N_2(v)-N_2(V)\r) = \eta\l\{{(\eta
    v)_\sigma}^2-{v_\sigma}^2\r\} + \tfrac12 \eta (\eta^2-1)v^2.
  \]
  Using the pointwise estimates for $v$ and $v_\sigma$ with the fact
  that $\eta\l (N_2(v)-N_2(V)\r)$ is supported in the region
  $e^{\delta\tau/2}\le |\sigma|\le \frac54 e^{\delta\tau/2}$, one
  arrives at the stated estimate for the $\fh$ norm of $\eta\l
  (N_2(v)-N_2(V)\r)$.
\end{proof}

We therefore have
\[
\|h_2\|^2_\fh a_2' (\tau) =\l (h_2, \eta N_2(V)\r)_\fh + \l (h_2,
\eta\l (N_2(v)-N_2(V)\r)\r)_\fh + \l (h_2, \eta N_3(v)+E\r )_\fh.
\]
If we write $N_2[v,w]=-v_\sigma w_\sigma-\frac12 vw$ for the natural
symmetric bilinear expression with $N_2(v)=N_2[v,v]$, then the first
term above can be written as
\begin{align*}
  \l (h_2, \eta N_2(V)\r)_\fh
  & = \l (\eta h_2, N_2 (a_2h_2 + a_0h_0 +W)\r)_\fh \\
  & = (a_2)^2 \l (\eta h_2, N_2(h_2)\r)_\fh \\
  & \quad +2a_2    \l (\eta h_2, N_2[h_2, a_0h_0+W]\r)_\fh \\
  & \quad +\l (\eta h_2, N_2(a_0h_0+W)\r)_\fh,
\end{align*}
which is a quadratic polynomial in $a_2$. Adding in the omitted error
terms, we reach the following observation.
\begin{lemma}\label{lem:a2-evolution}
  \[
  a_2'(\tau) = K(\tau) a_2(\tau)^2 + 2L(\tau) a_2(\tau) + M(\tau),
  \]
  in which
  \begin{align*}
    K(\tau) &= \frac1{\|h_2\| ^2}\l (\eta h_2, N_2 (h_2)\r)_\fh,\\
    L(\tau) &= \frac1{\|h_2\| ^2}\l (\eta h_2, N_2[h_2, a_0h_0+W]\r)_\fh, \\
    M(\tau) &= \frac1{\|h_2\| ^2} \l\{ \l (\eta h_2,
    N_2(a_0h_0+W)\r)_\fh +\l (h_2, \eta\l (N_2(v)-N_2(V)\r)\r)_\fh +
    \l (h_2, \eta N_3(v)+E\r )_\fh \r\}.
  \end{align*}
\end{lemma}

\begin{lemma}\label{lem:K-estimate}
  \[
  K(\tau) = -8 + \cO\l(e^{-e^{\delta\tau/2}}\r).
  \]
\end{lemma}
\begin{proof}
  Since $N_{2}(v)=-v_{\sigma}^{2}-\frac12 v^{2}$, the explicit
  expressions (\ref{eq:h024-explicit}) for $h_0$, $h_2$, and $h_4$
  give us
  \[
  N_2(h_2) = -\tfrac12\sigma^4-2\,\sigma^2-2 = -\tfrac12 h_4 -8\, h_2
  -12\, h_0.
  \]
  Hence $(h_2, N_2(h_2))_\fh = -8\|h_2\| ^2$. This implies
  \[
  K(\tau) = \frac{(h_2, N_2(h_2))_\fh}{\|h_2\| ^2}
  +\frac{((1-\eta)h_2, N_2(h_2))_\fh}{\|h_2\| ^2}.
  \]
  The first term is $-8$. The other term can be written as the
  integral of a function with fixed polynomial growth, which is
  supported in the region $|\sigma|\geq e^{\delta\tau/2}$. This leads
  quickly to the stated estimate.
\end{proof}
\begin{lemma}
  $L(\tau) = \cO\l(\tau^{-2}\r)$, i.e.\ there is a constant
  $C=C(F^*,n)$ such that for all $\tau\geq\tau_0+1$ one has
  $|\tau^2L(\tau)|\leq C$.
\end{lemma}
\begin{proof}
  Using (bi)linearity of $N_2[v,\tilde v]$, we get
  \begin{align*}
    \|h_2\|^2 L(\tau)
    &= \l(\eta h_2, N_2[h_2, a_0h_0+W]\r)_\fh \\
    &= a_0(\tau)\l(h_2,\eta N_2 (h_2,h_0)\r)_\fh +\l(h_2,\eta N_2
    (h_2, W)\r)_\fh.
  \end{align*}
  The first term is $\cO(\tau^{-2})$. For the other term, we have
  \begin{align*}
    \l| \l(h_2,\eta N_2 (h_2, W)\r)_\fh \r|
    & = \l|\l (\eta h_2, h_2'(\sigma)W_\sigma+\tfrac12 h_2 W \r)_\fh \r| \\
    & \le \|\eta h_2h_2'\|  \|W_\sigma\|  + C\|W\|  \\
    & \le \frac C{\tau^2}
  \end{align*}
  by our estimates for $\|W\| $ and $\|\pd_\sigma W\| $.
\end{proof}

\begin{lemma}\label{lem:M-bound}
  \[
  M(\tau) \leq C\left\{ \frac{e^{\tau_0-\tau}}{\tau}
    \|W(\cdot,\tau_0)\| +\frac1{\tau^3} \right\}
  \]
\end{lemma}
In fact one can show that $M(\tau) \leq C_\alpha\tau^{-(4-\alpha)}$
for any $\alpha>0$, but we will only need the estimate with
$\alpha=1$.
\begin{proof}
  In Lemma \ref{lem:a2-evolution}, $\|h_2\| ^2 M(\tau)$ is defined as
  the sum of three terms. We now estimate them one by one.

  The middle term is
  \[
  \l| (h_2, \eta (N_2(v)-N_2(V)))_\fh\r| \le \|h_2\| \|\eta
  (N_2(v)-N_2(V))\| \le C e^{-e^{\delta\tau/2}},
  \]
  by Lemma \ref{lem:N2-N2-estimate}.

  The third term is bounded by
  \[
  \l| (h_2, \eta N_3(v)+E)_\fh\r| \leq C\l (\|\eta N_3(v)\| + \|E\|
  \r) \leq \frac C{\tau^3},
  \]
  because of Lemma \ref{lem:N3-estimate} and the super-exponential
  estimate for $E$ in Lemma \ref{lem:E-estimate}.

  The first term can itself be split into three terms:
  \begin{equation}
    \label{eq:M-est-1}
    \bigl(\eta h_2, N_2 (a_0 h_2 +W)\bigr)_\fh
    = a_0(\tau)^2 \bigl(\eta h_2, N_2(h_2)\bigr)_\fh
    + 2a_0 (\tau) \bigl(\eta h_2, N_2[h_2, W]\bigr)_\fh
    + \bigl(\eta h_2, N_2(W)\bigr)_\fh .
  \end{equation}
  The first of these terms is $\cO (\tau^{-4})$, since $a_0(\tau)=\cO
  (\tau^{-2})$.

  From
  \begin{align*}
    \l (\eta h_2, N_2[h_2, W]\r)_\fh
    &=\l (\eta h_2, -h_2'(\sigma)W_\sigma -\tfrac12 h_2 W\r)_\fh \\
    &=-\bigl (\eta h_2h_2', W_\sigma\bigr)_\fh
    -\frac12\bigl (\eta h_2 (\sigma)^2, W)_\fh \\
    &=\cO\l (\|W_\sigma\|  + \|W\|  \r) \\
    &=\cO\l (\tau^{-2}\r)
  \end{align*}
  and $a_0 (\tau)=\cO(\tau^{-2})$, we deduce that the middle term in
  (\ref{eq:M-est-1}) is also $\cO (\tau^{-4})$.

  We are left with the last term in (\ref{eq:M-est-1}). This term can
  be expanded as
  \begin{equation}
    \label{eq:M-est-2}
    \bigl (\eta h_2, N_2(W)\bigr)_\fh = 
    -\bigl (\eta h_2, W_\sigma^2\bigr)_\fh
    -\bigl (\eta h_2, \tfrac12 W^2)_\fh.
  \end{equation}
  The second term satisfies
  \begin{align*}
    \l| \int_\R e^{-\sigma^2/4} \eta
    h_2(\sigma)W(\tau,\sigma)^2\,d\sigma \r|
    &\leq C \l (\|W\| ^2 + \|\sigma W\| ^2 \r ) \\
    &\leq C\l\{e^{2(\tau_0-\tau)} \|W(\cdot,\tau_0)\| ^2 + \frac
    1{\tau^4}\r\},
  \end{align*}
  by Lemmas \ref{lem:W-decay} and \ref{lem:sigmaW-bound}.

  We bound the first term in (\ref{eq:M-est-2}) as follows:
  \begin{align*}
    | \bigl (\eta h_2, W_\sigma^2\bigr)_\fh | & = \l|\int_\R
    e^{-\sigma^2/4} \eta h_2(\sigma)W_\sigma^2\,
    d\sigma\r| \\
    &\le C \l|\int_\R e^{-\sigma^2/4} (1+\sigma^2) W_\sigma^2\,
    d\sigma\r| \\
    &\le C \| (1+\sigma)^2 W_\sigma \| \|W_\sigma\| .
  \end{align*}
  The pointwise estimate (\ref{eq:Wsigma-bounds}) bounds the first
  factor by
  \[ \| (1+\sigma)^2 W_\sigma \| \leq \frac C\tau, \]
  while Lemma \ref{eq:W-deriv-est} implies
  \[
  \|W_\sigma\| \leq C\l\{e^{\tau_0-\tau} \|W(\cdot,\tau_0)\| + \frac
  1{\tau^2}\r\}.
  \]
  Combined, these last two estimates imply Lemma \ref{lem:M-bound}.
\end{proof}

\subsection{Dichotomy for the decay of $a_2(\tau)$}
We have shown that
\begin{equation}
  \label{eq:a2-evolution}
  \frac{da_2(\tau)}{d\tau}
  = -8 a_2(\tau)^2 + M^*(\tau),
\end{equation}
where
\[
\l|M^*(\tau)\r| \leq C(1+W_0)\tau^{-3}.
\]
Hence
\[
\tau\frac{d (\tau a_2(\tau))}{d\tau} = (\tau a_2 (\tau)) - 8 (\tau a_2
(\tau))^2 + \tau^2M^*(\tau).
\]
\begin{lemma}
  The coefficient $a_2(\tau)$ satifies either
  \begin{equation}
    \label{eq:a2-fastdecay}
    a_2 (\tau) \leq C(1+W_0)\tau^{-2}
  \end{equation}
  for all $\tau\geq\tau_0+1$, or else
  \begin{equation}
    \label{eq:a2-genericdecay}
    \lim_{\tau\to\infty} \tau a_2(\tau) = \frac18.
  \end{equation}
\end{lemma}

\subsection{Pinching time}

We will show below that alternative (\ref{eq:a2-genericdecay}) does
indeed occur, namely that there exists an open set of solutions
matching the asymptotic profile we have constructed. Our first step is
to improve the estimate obtained in Lemma 6.1 of \cite{AK}. Let
$r(t)\equiv\psi(0,t)$ denote the radius of the neck, and recall that
\[
K=\frac{\psi_{ss}}{\psi}\qquad\text{and}\qquad L=\frac{1-\psi_{s}^{2}}
{\psi^{2}}
\]
denote the sectional curvatures. Define $\lambda=2-\log L_{\min}(0)$,
so that we may write the cylindricality $F$ as
\[
F=\frac{K}{L}(\lambda+\log
L)=\frac{\psi\psi_{ss}}{1-\psi_{s}^{2}}\left(
  \lambda+\log\frac{1-\psi_{s}^{2}}{\psi^{2}}\right) .
\]
Recall that there exists $F^{\ast}\geq n-1$ depending only on $g(0)$
such that $0\leq F\leq F^{\ast}$ on the neck. Writing the evolution
equation for the radius
\[
\frac{dr}{dt}=\left[ \psi_{ss}-(n-1)\frac{1-\psi_{s}^{2}}{\psi}\right]
_{s=0}=\psi_{ss}-\frac{n-1}{r}
\]
in the form
\begin{equation}
  \frac{d}{dt}(r^{2})=-2(n-1)+2\frac{F}{\lambda+\log\frac{1}{r^{2}}}
\end{equation}
then proves that
\begin{equation}
  -2(n-1)\leq\frac{d}{dt}(r^{2})\leq-2(n-1)+2\frac{F^{\ast}}{\lambda+\log
    \frac{1}{r^{2}}}. \label{refined-shrink-estimate}
\end{equation}
This observation lets us prove that the singularity time depends
continuously on the initial conditions.

\begin{lemma}
  \label{PinchingContinuity}For rotationally symmetric metrics $g(t)$ that
  develop a neckpinch, the pinching time $T$ is a continuous function
  of the initial metric $g_{0}$.
\end{lemma}

\begin{proof}
  Given a metric $g_{0}$, let $T_{0}=T(g_{0})$ denote the singularity
  time of the solution $g(t)$ with initial data $g_{0}$.

  To prove lower semicontinuity, let $T_{-}<T_{0}$ be given. Then
  $\left\{ g(t):0\leq t\leq T_{-}\right\} $ is a smooth solution of
  Ricci flow. Because regular solutions of parabolic equations depend
  continuously on their initial data, a solution $\tilde{g}(t)$ will
  remain smooth for $0\leq t\leq T_{-}$, provided that its initial
  data $\tilde{g}_{0}$ is sufficiently close to $g_{0}$. This implies
  that the singularity time $\tilde{T}$ of $\tilde {g}(\cdot)$
  satisfies $\tilde{T}>T_{-}$.

  To prove upper semicontinuity, let $T_{+}>T_{0}$ be given. At any
  time $T_{\ast}<T_{0}$, observe that for all $\tilde{g}_{0}$
  sufficiently close to $g_{0}$, the radius $\tilde{r}$ of the
  solution $\tilde{g}(t)$ satisfies
  $\tilde{r}(T_{\ast})^{2}\leq2r(T_{\ast})^{2}$, and its
  cylindricality $\tilde{F}$ satisfies
  $\tilde{F}\leq\tilde{F}^{\ast}\leq2F^{\ast}$. By taking $T_{\ast}$
  close enough to $T_{0}$, one may by (\ref{refined-shrink-estimate})
  ensure that any solution $\tilde{g}(t)$ with initial data
  $\tilde{g}_{0}$ satisfies
  \[
  \frac{d}{dt}(\tilde{r}^{2})\leq-2(n-1)+\frac{4F^{\ast}}{\lambda+\log\frac
    {1}{2r(T_{\ast})^{2}}}\leq-(n-1)
  \]
  for all $t\in(T_{\ast},\tilde{T})$, where $\tilde{T}$ is the
  pinching time of $\tilde{g}(t)$. So for $t\in(T_{\ast},\tilde{T})$,
  it follows that
  \begin{align*}
    \tilde{r}(t)^{2}  &  \leq\tilde{r}(T_{\ast})^{2}-(n-1)(t-T_{\ast})\\
    &  \leq2r(T_{\ast})^{2}-(n-1)(t-T_{\ast})\\
    & \leq4(n-1)(T_{0}-T_{\ast})-(n-1)(t-T_{\ast}).
  \end{align*}
  Here we used the fact that $r(t)^{2}\leq2(n-1)(T_{0}-t)$. Choose
  $T_{\ast }<T_{0}$ closer to $T_{0}$ if necessary so that
  $T_{\ast}+4(T_{0}-T_{\ast })<T_{+}$. Then if $t\geq
  T_{\ast}+4(T_{0}-T_{\ast})$, one has $\tilde {r}(t)^{2}\leq0$. This
  implies that $\tilde{T}\leq T_{\ast}+4(T_{0}-T_{\ast })<T_{+}$.
\end{proof}

\subsection{Initial data}

Given a target time $T$, we now construct initial data $u=1+v$ at
$t=t_{0}<T$ that will become singular at $T$ and satisfy the
asymptotic profile (\ref{eq:a2-genericdecay}). In the next subsection,
we shall obtain appropriate estimates for the function
$V(\sigma,\tau)=\eta(e^{-\delta\tau /2}\sigma)\cdot
v(\sigma,\tau)=\sum_{k=0}^{\infty}a_{2k}(\tau)h_{2k}(\sigma)$
corresponding to this construction.

Let $\varepsilon\in(0,\frac{1}{4\sqrt{n-1}})$ be a constant to be
chosen later. Near the developing neckpinch, we impose the initial
profile
\[
u=1+\alpha_{0}h_{0}(\sigma)+\alpha_{2}h_{2}(\sigma),
\]
where $\alpha_{0}$ and $\alpha_{2}>0$ are constants to be chosen
later, and $\sigma(s,t)=s/\sqrt{T-t}$. That is, we prescribe
$\psi(s,t_{0})=R_{n} (t_{0})\cdot u(\sigma,t_{0})$ as
\begin{equation}
  \psi(s,t_{0})=\sqrt{2(n-1)}\sqrt{T-t_{0}}\left[  (1+\alpha_{0}-2\alpha
    _{2})+\alpha_{2}\frac{s^{2}}{T-t_{0}}\right]  \label{initial-profile}
\end{equation}
on the region $\mathcal{N}=\left\{ s:|s|\leq s_{0}\right\} $, where
$s_{0}$ is given by
\[
s_{0}=\varepsilon\sqrt{T-t_{0}}\log\frac{1}{T-t_{0}}.
\]
We will require
\begin{equation}
  -\varepsilon\leq\alpha_{0}\leq2\alpha_{2} \label{alpha_0-condition}
\end{equation}
and
\begin{equation}
  \frac{\varepsilon}{\tau_{0}}\leq\alpha_{2}\leq\frac{1}{\tau_{0}},
  \label{alpha_2-condition}
\end{equation}
where
\[
\tau_{0}=\log\frac{1}{T-t_{0}}\gg1.
\]

The radius of a neck of the form (\ref{initial-profile}) is
\begin{equation}
  r_{0}=\psi(0,t_{0})=\sqrt{2(n-1)}\sqrt{T-t_{0}}(1+\alpha_{0}-2\alpha_{2}).
  \label{initial-radius}
\end{equation}
Thus
\[
\psi=r_{0}+\sqrt{2(n-1)}\frac{\alpha_{2}s^{2}}{\sqrt{T-t_{0}}}.
\]
Since
\[
\psi_{s}=2\sqrt{2(n-1)}\frac{\alpha_{2}s}{\sqrt{T-t_{0}}},
\]
we have $|\psi_{s}|\leq\psi_{s}(s_{0},t_{0})$ on $\mathcal{N}$, where
$\psi_{s}(s_{0},t_{0})=2\varepsilon\sqrt{2(n-1)}\cdot\tau_{0}\alpha_{2}$
satisfies
\[
\varepsilon^{2}\leq\frac{\psi_{s}(s_{0},t_{0})}{2\sqrt{2(n-1)}}\leq
\varepsilon.
\]
Outside $\mathcal{N}$, we extend $\psi$ linearly with slope $\psi_{s}
(s_{0},t_{0})>0$ until we reach the height $\psi=1$, where we glue on
a standard profile that is smooth at the poles. As in \S 8 of
\cite{AK}, this can be done with $\left\vert \psi_{s}\right\vert
\leq1$ everywhere. Notice that the initial profile $\psi(\cdot,t_{0})$
is smooth except on $\partial\mathcal{N}$; the corresponding solution
of the Ricci flow will be smooth for all $t>t_{0}$ that it exists.

We next obtain a formula for the upper bound $F^{\ast}$ on the
cylindricality $F$ of a neck of the form (\ref{initial-profile}),
where $\alpha_{0}$ and $\alpha_{2}$ satisfy (\ref{alpha_0-condition})
and (\ref{alpha_2-condition}), respectively. Note that $F\leq0$ off
$\mathcal{N}$. On $\mathcal{N}$, we have
\begin{align*}
  \psi\leq\psi(s_{0},t_{0}) &
  =r_{0}+\sqrt{2(n-1)}\frac{\alpha_{2}s_{0}^{2}
  }{\sqrt{T-t_{0}}}\\
  & \leq\sqrt{2(n-1)}(1+\varepsilon^{2}\tau_{0})\sqrt{T-t_{0}}
\end{align*}
and
\[
\psi_{ss}=2\sqrt{2(n-1)}\frac{\alpha_{2}}{\sqrt{T-t_{0}}}\leq\frac
{2\sqrt{2(n-1)}}{\tau_{0}\sqrt{T-t_{0}}}.
\]
Independent of $\varepsilon\leq\frac{1}{4}$, we may choose $t_{\ast}$
close enough to $T$ so that for any $t_{0}\in\lbrack t_{\ast},T)$, the
radius is bounded from below by
\[
r_{0}\geq\sqrt{2(n-1)}\sqrt{T-t_{0}}\left(
  1-\varepsilon-\frac{2}{\tau_{0} }\right)
\geq\frac{\sqrt{2(n-1)}}{2}\sqrt{T-t_{0}},
\]
and so that the sectional curvature $L$ is bounded from below on
$\mathcal{N}$ by
\[
L=\frac{1-\psi_{s}^{2}}{\psi^{2}} \geq\frac{1-8(n-1)\varepsilon^{2}}
{\psi(s_{0},t_{0})^{2}}
\geq\frac{1}{4(n-1)(T-t_{0})(1+\varepsilon^2\tau_{0})^{2}} \geq e^{2},
\]
provided $\tau_0>\tau_0^*(n,\varepsilon)$.  Here we used the fact that
$\lim_{x\searrow0}(x\log\frac{1}{x})=0$. By Lemma 9.1 of \cite{AK},
the bound $L\geq e^{2}$ lets us assume $\lambda=0$.  Combining the
estimates above, we conclude that there exists $C=C(n,t_{\ast})$ such
that for any $t_{0}\in\lbrack t_{\ast},T)$, one has
\begin{align*}
  F=\frac{K}{L}\log L & =\frac{\psi\psi_{ss}}{1-\psi_{s}^{2}}\log\frac
  {1-\psi_{s}^{2}}{\psi^{2}}\\
  &  \leq\frac{2}{1-8(n-1)\varepsilon^{2}}\psi\psi_{ss}\log\frac{1}{r_{0}}\\
  &
  \leq16(n-1)\frac{1+\varepsilon^{2}\tau_{0}}{\tau_{0}}(\log2+\frac{\tau_{0}
  }{2})\\
  & \leq C\varepsilon^{2}\tau_{0}.
\end{align*}
Hence we may define
\begin{equation}
  F^{\ast}=(n-1)+C\varepsilon^{2}\tau_{0}. \label{bound-F-star}
\end{equation}
This formula for $F^{\ast}$ lets us adjust our constants in order to
force pinching at $T$.

\begin{lemma}
  There exist $\varepsilon_{\ast}$ and $t_{\ast}$ depending only on
  $n$ such that for any $\varepsilon\in(0,\varepsilon_{\ast}]$, all
  initial times $t_{0}\in\lbrack t_{\ast},T)$, and any $\alpha_{2}$
  satisfying (\ref{alpha_2-condition}), there exists some $\alpha_{0}$
  satisfying (\ref{alpha_0-condition}) such that the solution with
  initial profile (\ref{initial-profile}) becomes singular exactly at
  $T$.
\end{lemma}

\begin{proof}
  Proposition 5.2 of \cite{AK} implies that $r$ is monotonically
  decreasing. So for all $t\geq t_{0}$ that the solution exists, one
  has
  \[
  r(t)^{2}\leq r_{0}^{2}\leq2(n-1)(T-t_{0}).
  \]
  As above, we fix $t_{\ast}<T$ so that (\ref{bound-F-star}) holds for
  all $t_{0}\in\lbrack t_{\ast},T)$. Then by
  (\ref{refined-shrink-estimate}), there exists $c>0$ depending only
  on $n$ such that
  \[
  -2(n-1)\leq\frac{d}{dt}(r^{2})\leq-2(n-1)+2\frac{n-1+C\varepsilon^{2}\tau_{0}
  }{\tau_{0}-c}.
  \]
  There exists $C_{\ast}=C_{\ast}(n,t_{\ast})$ such that
  $\frac{n-1+C\varepsilon ^{2}\tau_{0}}{\tau_{0}-c}\leq
  C_{\ast}\varepsilon^{2}$. Hence we may fix
  $\varepsilon_{\ast}=\varepsilon_{\ast}(C_{\ast})$ such that for all
  $\varepsilon\in(0,\varepsilon_{\ast}]$, one has
  \[
  -2(n-1)\leq\frac{d}{dt}(r^{2})\leq-2(n-1)\left( 1-\varepsilon\right).
  \]
  By Lemma 7.1 of \cite{AK}, the solution becomes singular at
  $\tilde{T}$ if and only if $r(\tilde{T})=0$. It follows that
  \[
  \frac{r_{0}^{2}}{2(n-1)}\leq\tilde{T}-t_{0}\leq\frac{1}{1-\varepsilon}
  \cdot\frac{r_{0}^{2}}{2(n-1)}.
  \]
  By (\ref{initial-radius}), this implies that
  \[
  (1+\alpha_{0}-2\alpha_{2})^{2}\leq\frac{\tilde{T}-t_{0}}{T-t_{0}}\leq
  \frac{(1+\alpha_{0}-2\alpha_{2})^{2}}{1-\varepsilon}.
  \]
  Thus if $\alpha_{0}-2\alpha_{2}>0$, then pinching occurs at
  $\tilde{T}>T$. On the other hand, if
  $-1<\alpha_{0}-2\alpha_{2}<-\varepsilon$, then pinching occurs at
  $\tilde{T}<T$. By Lemma \ref{PinchingContinuity}, pinching will
  occurs exactly at $T$ for some
  $\alpha_{0}\in\lbrack2\alpha_{2}-\varepsilon ,2\alpha_{2}]$.
  Clearly, this $\alpha_{0}$ satisfies (\ref{alpha_0-condition}).
\end{proof}

\subsection{The dominant asymptotic profile}

For $\tau\geq\tau_{0}$, let $v(\sigma,\tau)$ denote the solution of
the Ricci flow corresponding to the initial data
(\ref{initial-profile}). Then
\[
v(\sigma,\tau_{0})=\alpha_{0}h_{0}(\sigma)+\alpha_{2}h_{2}(\sigma
)\qquad\text{for}\qquad|\sigma|\leq\sigma_{0}=\sigma(s_{0},t_{0}
)=\varepsilon\tau_{0}.
\]
For $\tau\geq\tau_{0}$, define $V(\sigma,\tau)=\eta(e^{-\delta\tau/2}
\sigma)v(\sigma,\tau)$, noting that $V(\cdot,\tau)\in\mathcal{H}$
admits the orthogonal decomposition
\begin{equation}
  V(\sigma,\tau)=a_{0}(\tau)h_{0}(\sigma)+a_{2}(\tau)h_{2}(\sigma)+W(\sigma
  ,\tau). \label{VerySpecialSolution}
\end{equation}

Intuitively, one expects that $a_{0}$ and $a_{2}$ will be close to
$\alpha _{0}$ and $\alpha_{2}$, respectively, at $\tau=\tau_{0}$ and
for a short time thereafter. We now make this expectation precise.
Define $U\in\mathcal{H}$ by
\[
U=V(\cdot,\tau_{0})-\alpha_{0}h_{0}-\alpha_{2}h_{2}.
\]
At $\tau=\tau_{0}$, one has $|V_{\sigma}|=|v_{\sigma}|=|\psi_{s}
|/\sqrt{2(n-1)}\leq\varepsilon$ for $|\sigma|\leq\sigma_{0}$ and
$|V_{\sigma }|\leq1$ everywhere. Because
$|\alpha_{0}|,\alpha_{2}\in\lbrack0,1]$, there are universal constants
$C$ (which may change from line to line) such that
\[
\begin{array}
  [c]{ll}
  U(\sigma)=0 & \text{for all }|\sigma|<\sigma_{0}\\
  |U(\sigma)|\leq C(1+\sigma^{2}) & \text{for all }|\sigma|\geq\sigma_{0},
\end{array}
\]
hence such that
\begin{align*}
  \left\Vert U\right\Vert ^{2} & \leq C\int_{|\sigma|\geq\sigma_{0}}
  (1+\sigma^{2})^{2}\,e^{-\sigma^{2}/4}\,d\sigma\\
  & \leq Ce^{-\sigma_{0}^{2}/5}=Ce^{-\varepsilon^{2}\tau_{0}^{2}/5}.
\end{align*}
Because $\left\Vert W\right\Vert ^{2}=\left\langle W,W\right\rangle
=\left\langle W,U\right\rangle \leq\left\Vert W\right\Vert \left\Vert
  U\right\Vert $ at $\tau=\tau_{0}$, it follows that the constant
$W_0$ introduced in Lemma \ref{lem:W-decay} satisfies
\begin{equation}
  W_{0}=e^{\tau_{0}}\left\Vert W(\cdot,\tau_{0})\right\Vert \leq e^{\tau_{0}
  }\left\Vert U\right\Vert \leq C\exp(\tau_{0}-\frac{\varepsilon^{2}}{10}
  \tau_{0}^{2}). \label{W_0-bound}
\end{equation}
Similarly, one has
\[
a_{2}(\tau_{0})=\frac{1}{\left\Vert h_{2}\right\Vert ^{2}}\left\langle
  h_{2},V(\cdot,\tau_{0})\right\rangle =\alpha_{2}+\frac{1}{\left\Vert
    h_{2}\right\Vert ^{2}}\left\langle h_{2},U\right\rangle ,
\]
and hence
\begin{equation}
  |a_{2}(\tau_{0})-\alpha_{2}|\leq Ce^{-\varepsilon^{2}\tau_{0}^{2}/10}.
  \label{Another-a_2-estimate}
\end{equation}

We are now ready to prove that alternative (\ref{eq:a2-genericdecay})
occurs.  By (\ref{eq:a2-evolution}), one has
\[
\frac{d}{d\tau}a_{2}(\tau)=-8a_{2}(\tau)^{2}+M^{\ast}(\tau),
\]
where $|M^{\ast}(\tau)|\leq C(1+W_{0})\tau^{-3}$. By
(\ref{W_0-bound}), we may make $\tau_{0}$ so large that
$C(1+W_{0})\leq C_{0}$. Then
\[
-8a_{2}(\tau)^{2}-\frac{C_{0}}{\tau^{3}}\leq\frac{d}{d\tau}a_{2}(\tau
)\leq-8a_{2}(\tau)^{2}+\frac{C_{0}}{\tau^{3}}.
\]
By (\ref{alpha_2-condition}) and (\ref{Another-a_2-estimate}), we may
make $\tau_{0}$ larger if necessary so that
\[
\sqrt{\frac{C_{0}}{8\tau_{0}^{3}}} \ll\frac{\varepsilon}{2\tau_{0}}
\leq a_{2}(\tau_{0}) \leq\frac{2}{\tau_{0}}.
\]
As long as $a_{2}(\tau)^2>\frac{C_{0}}{8}\tau^{-3}$, the function
$a_{2} (\tau)$ will be decreasing. For such times, one has
\[
\frac{d}{d\tau}a_{2}(\tau)
\geq-8a_{2}(\tau_{0})^{2}-\frac{C_{0}}{\tau^{3}}
\geq-\frac{32}{\tau_{0}^{2}}-\frac{C_{0}}{\tau_{0}^{3}}
\geq-\frac{64}{\tau_{0}^{2}}
\]
provided that $\tau_{0}$ is large enough, and hence
\[
a_{2}(\tau)\geq a_{2}(\tau_{0})-\frac{64}{\tau_{0}^{2}}(\tau-\tau_{0}
)\geq\frac{\varepsilon}{2\tau_{0}}-\frac{64}{\tau_{0}^{2}}(\tau-\tau_{0}).
\]
In particular, given any $C>0$, one can make $\tau_{0}$ so large that
\[
a_{2}(\tau_{0}+1)
\geq\frac{\varepsilon}{2\tau_{0}}-\frac{64}{\tau_{0}^{2}} \geq
C\frac{1+W_{0}}{(\tau_{0}+1)^{2}}. 
\]
This justifies our assumption that $a_{2}$ was decreasing for
$\tau_{0} \leq\tau\leq\tau_{0}+1$. Moreover, if we choose $C$ as in
(\ref{eq:a2-fastdecay}) then it violates alternative
(\ref{eq:a2-fastdecay}). Hence we have proved the following result.

\begin{proposition}
  For any final time $T$, all sufficiently small $\varepsilon>0$, all
  initial times $t_{0}$ sufficiently close to $T$, and any
  $\alpha_{2}$ satisfying (\ref{alpha_2-condition}), there exists
  $\alpha_{0}$ satisfying (\ref{alpha_0-condition}) such that the
  solution (\ref{VerySpecialSolution}) with initial profile
  (\ref{initial-profile}) becomes singular at $T$ and satisfies the
  asymptotic profile
  \begin{equation}\label{eq:var-a2-genericdecay}
    a_{2}(\tau)=\frac{1}{8\tau}+o(\frac{1}{\tau})\qquad\text{as}\qquad
    \tau\rightarrow\infty.
  \end{equation}

\end{proposition}

\subsection{The intermediate region}

Consider a solution that forms a singularity at time $t=T$. Define $u$, $v$, etc.~as
above and assume that our solution satisfies (\ref{eq:var-a2-genericdecay}).

We have shown that
\begin{equation}
  v(\sigma,\tau)=\frac{\sigma^{2}-2}{8\tau}+o(\frac{1}{\tau})\qquad
  (\tau\rightarrow\infty), \label{eq:v-inner-expansion}%
\end{equation}
uniformly on any bounded interval $|\sigma|\leq A$.

Now consider the quantity
\[
q=u^{2}-1.
\]

\begin{lemma}\label{lem:intermediate-convergence}
  For any constants $\varepsilon>0$ and $B<\infty$, there exist $A$ and
  $\bar{\tau}$ such that
  \[
  (1-\varepsilon)\frac{\sigma^{2}}{4\tau} <
  q(\sigma, \tau)<(1+\varepsilon)\frac{\sigma^{2}}{4\tau}
  \]
  whenever $A\leq\sigma\leq B\sqrt\tau$ and times $\tau\geq \bar{\tau}$.
\end{lemma}

In other words, for $1\ll\sigma\leq B \sqrt\tau $ one has
\[
u = \sqrt{1+(1+o(1))\frac{\sigma^2}{4\tau}}
\approx\sqrt{1+\frac{\sigma^2}{4\tau}}.
\]

We shall prove this Lemma by computing the evolution equation for $q$ in the
variables $\tau$ and $\rho=\sigma/\surd\tau$, and finding suitable sub- and
supersolutions. The evolution equation (\ref{eq:q-rho-const}) turns out to be a small
perturbation of a first order linear equation, which is easily solved by following
its characteristics.

\subsubsection{Equation for $q$}
Since $q_{\sigma}=2uu_{\sigma}$ and $q_{\sigma\sigma}=2uu_{\sigma\sigma} +
2u_{\sigma}^{2}$, one has
\begin{equation}
  \label{eq:q-x-const}
  q_{\tau}|_{x}=q_{\sigma\sigma}+q+\frac{n-2}{2}\frac{q_{\sigma}^{2}}{q+1}.
\end{equation}
Here the time derivative is with constant $x$. Using Lemma \ref{lem:J}, we get
\begin{equation}\label{eq:q-s-const}
  \left.  q_{\tau}\right\vert_{\sigma}
  =q_{\sigma\sigma}-\frac{\sigma}{2}
  q_{\sigma}+q-nJ(\sigma,\tau)q_{\sigma}+\frac{n-2}{2}\frac{q_{\sigma}^{2}}{q+1},
\end{equation}
where $\left. q_{\tau}\right\vert _{\sigma}$ is the derivative of $q$ with $\sigma$
held constant and where
\[
J=\frac{u_{\sigma}}{u}+\int_{0}^{\sigma}\frac{u_{\sigma}^{2}}{u^{2}}d\sigma.
\]
We now regard $q$ as a function of the new variable $\rho=\sigma/\sqrt{\tau}$.
Because
\[
\left.  q_{\tau}\right\vert _{\sigma}=\left.  q_{\tau}\right\vert _{\rho
}-\frac{1}{2}\frac{\rho}{\tau}q_{\rho},
\]
one has
\begin{equation}
  \left.
    q_{\tau}\right\vert _{\rho}-\frac{1}{\tau}q_{\rho\rho}+\frac{\rho}{2}q_{\rho}-q
  =\frac{1}{\tau}\left\{
    \frac{\rho}{2}q_{\rho}
    -n\sqrt{\tau}Jq_{\rho}
    +\frac{n-2}{2}\frac{q_{\rho}^{2}}{q+1}
  \right\}  .
\label{eq:q-rho-const}
\end{equation}

\subsubsection{A subsolution for (\ref{eq:q-rho-const})}

We begin with some simple estimates. Fix some $A_{0}\gg1$ and recall that $\tau_{0}$
is the smallest time under consideration in this paper. Given $B>1$, we shall choose
$A\geq A_{0}$ and $\tau_{1}\geq\tau_{0}$ large enough to make everything work. The
region $A/\sqrt{\tau}\leq\rho\leq B^{2}$ corresponds to $A\leq\sigma\leq
B^{2}\sqrt{\tau}$. By Lemma \ref{lem:pointwise-u-estimates}, there exist
$C_{0}=C_{0}(B)$ and $C_{1} =C_{1}(A_{0})$ such that for all $A\leq\sigma\leq
B^{2}\sqrt{\tau}$, one has
\[
u\leq1+C_{0}\frac{\sigma^{2}}{\tau}=1+C_{0}\rho^{2}\leq1+C_{0}B^{4}
\]
and
\[
0\leq u_{\sigma}\leq C_{0}\frac{1+\sigma}{\tau}\leq C_{1}\frac{\sigma}{\tau
}=C_{1}\frac{\rho}{\sqrt{\tau}}.
\]
Thus we have
\[
0\leq u_{\rho}\leq C_{1}\rho
\]
and
\[
0\leq\rho q_{\rho}=2\rho u_{\rho}u\leq C(A_{0},B)
\]
for $A/\sqrt{\tau}\leq\rho\leq B^{2}$. Since $u\geq1-C/\tau$ for $C=C(n,F^{\ast})$,
we may assume without loss of generality that $\tau_{0}\gg1$ is so large that
\[
\frac{q_{\rho}^{2}}{q+1}=\frac{q_{\rho}^{2}}{u^{2}}\leq2q_{\rho}^{2}\leq
C(A_{0},B)
\]
for $A/\sqrt{\tau}\leq\rho\leq B^{2}$ and $\tau\geq\tau_{0}$. Since
$n\sqrt{\tau}Jq_{\rho}\geq0$ for all $0\leq\sigma\leq e^{\delta\tau}$, it follows
that
\begin{equation}
\frac{\rho}{2}q_{\rho}-n\sqrt{\tau}Jq_{\rho}+\frac{n-2}{2}\frac{q_{\rho}^{2}
}{q+1}\leq C(A_{0},B).\label{LqEstimate}
\end{equation}

We introduce the linear differential operator
\[
\mathcal{L}[f]
\doteqdot
f_{\tau}|_{\rho}-\frac{1}{\tau}f_{\rho\rho}+\frac{\rho}{2}f_{\rho}-f.
\]
Then (\ref{eq:q-rho-const}) and (\ref{LqEstimate}) imply that $q=u^{2}-1$ satisfies
\[
|\mathcal{L}[q]|\leq\frac{C(A_{0},B)}{\tau}.
\]
whenever\ $A/\sqrt{\tau}\leq\rho\leq B^{2}$ and $\tau\geq\tau_{0}$.

Let $Q$ be a smooth function with support in $[0,1]$ such that $Q(r)\leq r^{2}/4$,
and define
\[
K_{0}=\max_{0\leq r\leq1}|Q^{\prime\prime}(r)|.
\]
For $\tau_{1}\geq\tau_{0}$, define
\begin{equation}\label{eq:q0-def}
  q_{0}(\rho,\tau)
  =\rho_1(\tau)^2 Q\bigl(\rho/\rho_1(\tau)\bigr),
\end{equation}
with 
\[
\rho_1(\tau) \doteqdot e^{(\tau-\tau_1)/2}\frac A{\sqrt\tau_1}.
\]
This function is almost in the nullspace of the linear operator $\mathcal{L}$;
indeed, a straightforward computation shows that
\[
\mathcal{L}[q_{0}]
=-\frac{1}{\tau}\frac{\partial^{2}q_{0}}{\partial\rho^{2}}
=-\frac{1}{\tau}Q^{\prime\prime}\bigl(\rho/\rho_1(\tau)\bigr),
\]
so that $|\mathcal{L}[q_{0}]|\leq K_{0}/\tau$.

Now for $K\geq K_{0}$ to be determined and any $\theta\in(1/2,1)$, consider
\begin{equation}
  q_{-}(\rho,\tau)=\theta q_{0}(\rho,\tau)+\frac{K}{\tau}.
\end{equation}
One then has
\[
\mathcal{L}[q_{-}]
=\theta\mathcal{L}[q_{0}]+K\mathcal{L}[\tau^{-1}]
<\theta\frac{K_{0}}{\tau}-\frac{K}{\tau}<\frac{K_{0}-K}{\tau}\leq0,
\]
so that $q_{-}$ is a subsolution for all $\tau\geq\tau_{0}$. Moreover, if we choose
$K\geq K_{0}+C(A_{0},B)$, then
\begin{equation}\label{SubsolutionEstimate}
  \mathcal{L}[q_{-}-q]\leq\frac{C(A_{0},B)+K_{0}-K}{\tau}
  \leq0
\end{equation}
for $A/\sqrt{\tau}\leq\rho\leq B^{2}$ and $\tau\geq\tau_{0}$.

Define the region
\begin{equation}\label{eq:Omega-def}
  \Omega=\left\{  (\rho,\tau) :
    \frac{A}{\sqrt{\tau}}\leq \rho \leq \rho_1(\tau),
    \;\tau_{1}\leq\tau\leq\omega
  \right\}  ,
\end{equation}
where $\omega$ is determined by $\rho_1(\omega) = B^2$, i.e.
\[
\omega=\tau_{1}+\log\left(  B^{4}\tau_{1}/A^{2}\right)  .
\]
We will choose $A$ and $\tau_{1}$ so that $B^{4}\tau_{1}>A^{2}$. First, we compare
the values of $q$ and $q_{-}$ along the parabolic boundary of $\Omega$.

Along the left edge $\rho=A/\sqrt{\tau}$, $\tau\geq\tau_{1}$, it follows from
(\ref{eq:v-inner-expansion}) that
\begin{equation}
  q(A/\sqrt{\tau},\tau)=u^{2}(A,\tau)-1=\frac{A^{2}-2}{4\tau}+o(\frac{1}{\tau
  })\geq(1-\varepsilon_{0})\frac{A^{2}}{4\tau},\label{qFromBelow}
\end{equation}
where $\varepsilon_{0}=\varepsilon_{0}(A_{0},\tau_{1})$. On the other hand,
there exists $\varepsilon_{1}=\varepsilon_{1}(A,K,\tau_{0})$ such that
\[
q_{-}(\rho,\tau)\leq\frac{\theta}{4}\rho^{2}+\frac{K}{\tau}\leq(1+\varepsilon
_{1})\frac{\theta}{4}\rho^{2}
\]
for $\rho\geq A/\sqrt{\tau}$. So $q_{-}<q$ at $(A/\sqrt{\tau},\tau)$ provided
that
\[
\frac{1}{2}\leq\theta<\frac{1-\varepsilon_{0}}{1+\varepsilon_{1}}.
\]

Along the right edge $\rho=\rho_1(\tau)$, one has
\[
q_{-}(\rho_1(\tau),\tau)=\rho_1(\tau)^2 Q(1) + \frac{K}{\tau}=\frac{K}{\tau},
\]
because $Q$ is supported in $[0,1]$. But since $q_{\rho}\geq0$, estimate
(\ref{qFromBelow}) implies that
\[
q(\rho_1(\tau),\tau)
 = q(e^{(\tau-\tau_{1})/2}A/\sqrt{\tau_{1}},\tau)
 \geq q(A/\sqrt\tau,\tau)
 \geq \frac{1-\varepsilon_{0}}{4}\frac{A^{2}}{\tau},
\]
where $\varepsilon_{0}=\varepsilon_{0}(A_{0},\tau_{1})$. So $q_{-}\leq q$ at
$(\rho_1(\tau), \tau)$ provided that
\[
 A\geq2\sqrt{\frac{K}{1-\varepsilon_{0}}}.
\]

The bottom of $\Omega$ is the single point $(A/\sqrt{\tau_{1}},\tau_{1})$. So the
choice $\theta<(1-\varepsilon_{0})/(1+\varepsilon_{1})$ ensures that $q_{-}<q$ there
as well.

Now choose $K$ depending on $A_{0}$, $B$, and $K_{0}$ to satisfy $K\geq
K_{0}+C(A_{0},B)$. Take $\tau_{1}\geq\tau_{0}$ large enough so that
$\varepsilon_{0}<1/2$ and choose $A\geq A_{0}$ depending on $K$ so that
$A\geq\sqrt{8K}$ and so that $\varepsilon_{1}>0$ is as small as desired. Then
increase $\tau_{1}$ if necessary so that $\tau_{1}\geq A^{2}/B^{4}$, making
$\varepsilon_{0}$ as small as desired. These choices make it possible to take
$\theta\in(1/2,1)$ as close to $1$ as one wishes and still satisfy
$\theta<(1-\varepsilon_{0})/(1+\varepsilon_{1})$. Furthermore, they ensure that
$q_{-}\leq q$ along the parabolic boundary of $\Omega$. Because
$e^{(\tau-\tau_{1})/2}A/\sqrt{\tau_{1}}\leq B^{2}$ for all $\tau\leq\omega$, it
follows from (\ref{SubsolutionEstimate}) that $q\geq q_{-}$ throughout $\Omega$.
Since for any $\varepsilon$, one may take $Q(r)=r^{2}/4$ for
$r\in(\varepsilon,1-\varepsilon)$, this implies in particular that
\[
q(\rho,\tau)>\theta\frac{\rho^{2}}{4}
\]
for all points $\rho$ such that $A/\sqrt{\tau}\leq\rho\leq B$ and times $\tau$ such
that
\[
2\log\left(  B\sqrt{\tau_{1}}/A\right)  \leq\tau-\tau_{1}\leq2[\log\left(
B\sqrt{\tau_{1}}/A\right)  +\log B].
\]
Our estimates only improve if we increase $\tau_{1}$ (without altering any other
choices); in particular, $\theta$ does not decrease. Since $\tau = \tau_{1} +
o(\tau_{1})$ as $\tau_{1}\rightarrow\infty$, we have proved the lower estimate in
Lemma \ref{lem:intermediate-convergence}.

\subsubsection{A supersolution}
A slight modification of the preceding construction will give an upper barrier for
$q$ on the same domain $\Omega$. One begins with the \textit{a priori} bound
\begin{equation}
  \label{eq:u-upperbound}
  u\leq C (1+\sigma^2/\tau)= C(1+\rho^2)
  \text{ and hence }
  q\leq C (1+\rho^2)\rho^2 = C_0(B)\rho^2,
\end{equation}
which holds for $A_0/\sqrt\tau\leq \rho\leq B^2$.

Choose a smooth function $\hq$ which satisfies $\hq(r)\geq r^2/4$, and define $\hat
K_0=\max_{[0,1]}|\hq''|$. Instead of requiring $\hq$ to be supported in $[0,1]$, we
impose
\[
\hq(1)=2C_0(B)
\text{ with $C_0(B)$ as in (\ref{eq:u-upperbound}). }
\]
For $\tau\geq \tau_1$, one introduces \( \hat q_{0}(\rho,\tau) =\rho_1(\tau)^2 \hat
Q\bigl(\rho/\rho_1(\tau)\bigr), \) as in (\ref{eq:q0-def}), and finds that
$|\mathcal{L}[\hat q_0]|\leq \hat K_0/\tau$. For $\hat K\geq\hat K_0$ and any
$\theta\in(1,\frac32)$, define
\[
 q_+(\rho,\tau) = \theta\hat q_0(\rho,\tau) - \hat K / \tau.
\]
Choosing $\hat K\geq \hat K_0+C(A_0,B)$, one finds as in (\ref{SubsolutionEstimate})
that $\cL[q_+-q]\geq0$ for $\tau\geq\tau_0$ and $A/\sqrt\tau\leq\rho\leq B^2$. Thus
$q_+$ is a supersolution, provided $q_+\geq q$ on the parabolic boundary of $\Omega$.

On the left edge of $\Omega$, one has, as in (\ref{qFromBelow}), \( q\leq
(1+\varepsilon_0)A^2/4\tau \). One also has \( q_+\geq (1-\varepsilon_1)\theta
A^2/4\tau \), where $\varepsilon_1=\varepsilon_1(A,K)$.

Hence one will have $q\leq q_+$ on the left edge if $\theta$ satisfies
\[
\frac{1+\varepsilon_0} {1-\varepsilon_1} <\theta < \frac 32.
\]

On the right edge of $\Omega$ one has $\rho=\rho_1(\tau)$, and thus
\[
q_+(\rho,\tau)
= \theta\rho_1(\tau)^2\hat Q(1) - \hat K/\tau
\geq \theta\bigl(1-\hat K/A^2\bigr) Q(1)\rho_1(\tau)^2
Q(1)\rho_1(\tau)^2,
\]
provided $A$ is chosen so large that $\theta(1-\hat K/A^2)>1$. On the
other hand we know that $q\leq C_0(B)\rho_1(\tau)^2$ holds on the
right edge of $\Omega$. Hence our condition $\hat Q(1)>C_0(B)$ implies
that $q_+>q$ on the right edge.

It now follows that $q\leq q_+$ on $\Omega$, and as with the subsolution,
this implies that $q(\rho,\tau) <\theta\rho^2/4$ for all
$\rho\in [A/\sqrt\tau, B]$ and $\tau\geq\tau_1+o(\tau_1)$.

\subsection{Convergence in the outer region}
\label{sec:outer}

Let $B$ and $\varepsilon>0$ be given constants.  Define
$A=A(\varepsilon,B)$ and $\bar\tau=\bar\tau(\varepsilon,B)$ as in
Lemma \ref{lem:intermediate-convergence}.  Then we have shown that
when
\begin{equation}
  \label{eq:B-curve}
  \sigma=B\sqrt\tau, \text{ i.e.~when } s^2 = \frac{B^2}4(T-t)\log\frac1{T-t},
\end{equation}
one has
\[
\sqrt{1+ (1-\varepsilon)B^2/4} \leq u\leq \sqrt{1+
  (1+\varepsilon)B^2/4}.
\]
For $\psi = R_n(t) u$, this therefore means that
\begin{equation}
  \label{eq:psi-on-outer-edge}
  \sqrt{T-t}\sqrt{1+ (1-\varepsilon)B^2/4}\leq \frac{\psi}{\sqrt{2(n-1)}}\leq
  \sqrt{T-t}\sqrt{1+ (1+\varepsilon)B^2/4}
\end{equation}
for all $t\in[\bar t, T)$, where $\bar t = T-e^{-\bar \tau}$.

We define the \emph{outer region} (with parameter $B$) to be the
portion of space-time given by
\[
\Gamma_B = \bigl\{(x, t) :\bar\tau\leq \tau<T, s(x,t)^2\geq
\frac{B^2}4(T-t)\log\frac1{T-t}\bigr\}.
\]
The boundary $\partial\Gamma_B$ of $\Gamma_B$ is given by
(\ref{eq:B-curve}).

For any given $t_{\ast}\in\lbrack\bar{t},T)$, we let $x_{\ast}$ be the $x$
coordinate of the point on $\partial\Gamma_{B}$ with $t=t_{\ast}$. We will now
show that $s(x_{\ast},t)$ and $\psi(x_{\ast},t)$ change very little for
$t_{\ast}\leq t<T$. This will directly imply the asymptotic description of
$\psi(s,T)$ for $s\rightarrow0$ in the introduction. Fix some $B_{0}>1$. We
begin with a simple uniform estimate for $\psi_{s}$ in the region $B_{0}%
\sqrt{\tau}\leq\sigma\leq\eta e^{\delta\tau}$.

\begin{lemma}
\label{lem:intermediate-usig-bound}For any $B_{0}>1$, there exist constants
$\gamma,\delta>0$ and $C<\infty$ depending only on $B_{0}$ and the initial data
such that%
\[
0\leq\psi_{s}(s,t)\leq\frac{C}{\sqrt{\tau}}%
\]
for all $\ \tau\geq\bar{\tau}(B_{0})$ and $B_{0}\sqrt{\tau}\leq\sigma
(s)\leq\gamma e^{\delta\tau}$.
\end{lemma}

\begin{proof}
Define $r(t)=\psi(0,t)$. By (\ref{eq:var-a2-genericdecay}), we have%
\[
r=\sqrt{2(n-1)(T-t)}\left\{  1-\frac{1}{4\tau}+o(\frac{1}{\tau})\right\}  .
\]
By (\ref{eq:psi-on-outer-edge}), we have%
\[
\psi(s,t)\geq\psi(\bar{s},t)\geq\sqrt{2(n-1)(T-t)}\sqrt{1-B_{0}^{2}/8}%
\]
for all $s\geq\bar{s}(B_{0})=B_{0}\tau\sqrt{T-t}$ and $\tau\geq\bar{\tau
}(B_{0})\geq\tau_{0}$. Now restricting our attention to the region of the neck
where $\psi\leq1/B_{0}$ and $0\leq\psi_{s}\leq1/B_{0}$, we may apply estimate
(34) of \cite{AK}, obtaining $C$ depending on $B_{0}$ and $F^{\ast}$ such that%
\[
\psi_{s}^{2}\leq C\log\left(  \frac{\log r}{\log\psi}\right).
\]
This estimate holds as long as $\sigma\leq\gamma e^{\delta\tau}$, i.e.~as long as $s\leq
\gamma(\sqrt{T-t})^{1-2\delta}\leq r^{1-2\delta}$. Let $c_{1}=\log
\sqrt{1-B_{0}^{2}/8}$. Then%
\[
\frac{\log r}{\log\psi}\leq\frac{\log\sqrt{2(n-1)}-\frac{\tau}{2}-\frac
{1}{4\tau}+o(\frac{1}{\tau})}{\log\sqrt{2(n-1)}-\frac{\tau}{2}+c_{1}}%
=1+\frac{2c_{1}}{\tau}+o(\frac{1}{\tau}).
\]
Hence $\tau\psi_{s}^{2}\leq C\tau\log\left(  \frac{\log r}{\log\psi}\right)
=C\left(  2c_{1}+o(1)\right)  $.
\end{proof}

We are now ready to describe $s(x_*, t)$ and $\psi(x_*, t)$.

\begin{lemma}
  For any $x_*$ the limits
  \[
  \lim_{t\to T}s(x_*, t) \text{ and }\lim_{t\to T}\psi(x_*, t) 
  \]
  exist.  Furthermore, a constant $C<\infty$ which does not depend on
  $B \geq B_0$ or $x_*$ exists for which one has
  \begin{align}
    \label{eq:s-variation}
    |s(x_*, t) - s(x_*,t_*)|
    &\leq C\bigl\{B^{-3}+|\log(T-t_*)|^{-3/2}\bigr\}\; s(x_*,t_*) \\
    \intertext{and}
    \label{eq:psi-variation}
    |\psi(x_*, t) - \psi(x_*,t_*)| & \leq \frac C{B^2} \psi(x_*,t_*)
  \end{align}
  for all $t\in[t_*, T)$.
\end{lemma}
\begin{proof}
  To estimate the change in $s(x_{*}, t)$ we use
   \begin{equation}
    \label{eq:outer-s-growth}
    \left.\frac{\partial s}{\partial t}\right|_{x=x_*} =
    n\frac{\psi_{s}(x_*, t)}{\psi(x_*, t)} +
    n\int_{0}^{s(x_*, t)}\left(\frac{\psi_s}{\psi}\right)^{2}ds,
  \end{equation}
  as in the proof of Lemma~\ref{TheNapkinLemma}.
  Let $t'\in(t_*,T)$ be given, and let $x'$ be such that $s(x',t')$ lies on
  $\partial\Gamma_B$.

 \psfrag{s}{$s$} 
 \psfrag{T}{$T$} 
 \psfrag{t}{$t$} 
 \psfrag{ts}{$t_*$} 
 \psfrag{tp}{$t'$}
 \psfrag{sxptp}{$s(x',t')$} 
 \psfrag{sxstp}{$s(x_*,t')$} 
 \psfrag{sxsts}{$s(x_*, t_*)$}
 \psfrag{bgb}{$\partial\Gamma_B$}
 \begin{center}
    \includegraphics{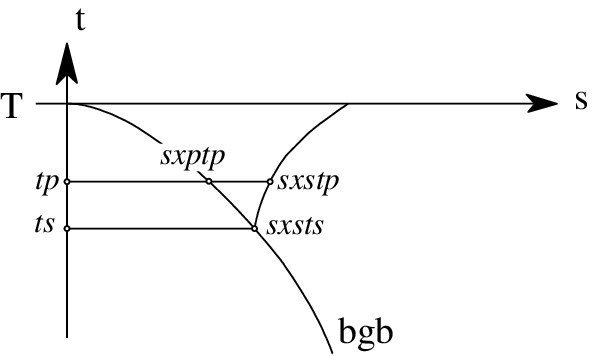}
 \end{center}

  Then at time $t'$, Lemma \ref{lem:intermediate-usig-bound} implies that
  \begin{equation}\label{eq:outer-psi-s-bound}
    0\leq\psi_s(s,t') \leq \frac C{\sqrt{\log \frac1{T-t'}}} \quad
    \text{for all }s\leq s'.
  \end{equation}
  Since $x_*\geq x'$, we then get
  \[
  \int_0^{s(x_*,t')}\left(\frac{\psi_{s}}{\psi}\right)^{2}ds
  =\int_0^{s(x',t')}\left(\frac{\psi_{s}}{\psi}\right)^{2}ds
  +\int_{s(x',t')}^{s(x_*,t')}\left(\frac{\psi_{s}}{\psi}\right)^{2}ds.
  \]
  In the first term, we use \eqref{eq:outer-psi-s-bound} to get
  \begin{align*}
    \int_0^{s(x',t')}\left(\frac{\psi_{s}}{\psi}\right)^{2}ds &\leq
    \frac C{\sqrt{\ln1/(T-t_*)}}\int_0^{s(x',t')}\frac{\psi_s}{\psi^2}ds \\
    &\leq \frac C{\sqrt{\ln1/(T-t_*)}}\frac1{\psi(0,t')} \\
    &\leq \frac C{\sqrt{(T-t')\ln1/(T-t_*)}}.
  \end{align*}
  In the second term, we have $x'<x<x_*$ so that
  \[
  \psi\geq\psi(x', t')\geq\sqrt{2(n-1)(T-t')\bigl(1+(1-\varepsilon)B^2/4\bigr)} \geq
  CB\sqrt{T-t'},
  \]
  with $C$ independent of $t$ or $B \geq B_0$.  On the other hand, we also have
  $\psi_s\leq1$, so that
  \begin{align*}
    \int_{s(x',t')}^{s(x_*,t')}\left(\frac{\psi_{s}}{\psi}\right)^{2}ds
    &\leq \int_{s(x',t')}^{s(x_*,t')}\left(\frac{\psi_{s}}{\psi^2}\right)ds \\
    &\leq \frac1{\psi(x',t')}\\
    &\leq \frac C{B\sqrt{T-t'}}.
  \end{align*}
  For the first term in (\ref{eq:outer-s-growth}), one has the same
  estimate, so that at any time $t'\in(t_*, T]$, one has
  \[
  0\leq \partial_ts(x_*, t') \leq C \bigl\{B^{-1}+|\log(T-t_*)|^{-1/2}\bigr\}
  (T-t')^{-1/2}.
  \]
  This is integrable in $t'\in(t_*, T)$, so that $\lim_{t\to T}s(x_*,
  t)$ must exist.  Integration with $t_*<t'<t$ then leads to
  \[
  |s(x_*,t)-s(x_*, t_*)|\leq C \bigl\{B^{-1}+|\log(T-t_*)|^{-1/2}\bigr\}\sqrt{T-t_*}.
  \]
  Recall (\ref{eq:B-curve}), apply $a^2b, ab^2\leq a^3+b^3$ to
  $a=B^{-1}$ and $b=|\log(T-t_*)|^{-1/2}$, and conclude
  \begin{align*}
     |s(x_*,t)-s(x_*, t_*)|
      &\leq C \bigl\{B^{-1}+|\log(T-t_*)|^{-1/2}\bigr\}
      \frac{s(x_*, t_*)}{B \sqrt{\ln1/(T-t_*)}} \\
      &\leq C \bigl\{B^{-3}+|\log(T-t_*)|^{-3/2}\bigr\}s(x_*, t_*),
  \end{align*}
  which proves the estimate (\ref{eq:s-variation}) for $s$.

  To establish the other estimate, we recall from Proposition 5.3 in
  \cite{AK} that there is a constant $C$ which only depends on the
  initial data and the dimension $n$ such that
  \[
  |2\psi\psi_t|\leq C
  \]
  always holds.  In particular, the constant $C$ does not depend on
  our choice of $B$.  In~\cite{AK}, we used this estimate to show that
  $\lim_{t\to T}\psi(x_*,t)$ exists. It also implies that
  \begin{equation}
    \label{eq:outer11}
      |\psi(x_*, t)^2-\psi(x_*,t_*)^2| \leq C(t-t_*)
  \end{equation}
  and hence, using (\ref{eq:psi-on-outer-edge}), that
  \[
  |\psi(x_*, t)-\psi(x_*,t_*)|
     \leq \frac{C(t-t_*)}{\psi(x_*, t)+\psi(x_*,t_*)}
     \leq \frac CB \sqrt{T-t_*}
     \leq \frac C{B^2}\psi(x_*, t_*)
  \]
  for all $t\in(t_*, T)$.
\end{proof}

\subsection{Asymptotics in the outer region}
Using (\ref{eq:psi-on-outer-edge}) and (\ref{eq:psi-variation}), we compute that
\begin{align*}
  \psi(x_*, t)
    &= \bigl(1+\cO(B^{-2})\bigr) \psi(x_*,t_*)\\
    &= \bigl(1+\cO(B^{-2})\bigr)
       \sqrt{2(n-1)\bigl[1+(1+\theta\varepsilon)B^2/4\bigr](T-t_*)}\\
    &= \bigl(1+\cO(B^{-2}+\varepsilon)\bigr)
       \sqrt{\frac{n-1}2}B\sqrt{T-t_*}
\end{align*}
for some $|\theta|\leq1$.

From the definition (\ref{eq:B-curve}), we see that $s=s(x_*, t_*)$
satisfies
\[
B\sqrt{T-t} = \frac s{\sqrt{-\log(T-t)}}.
\]
Hence, taking logarithms and abbreviating $\omega=T-t$,
\[
\log B +\tfrac12\ln\omega = \log s -\tfrac12\ln|\log\omega|,
\]
so that
\begin{align*}
  \log\omega
  &=2\ln s-\log|\log\omega|-2\ln B =\cO(\log s)\\
  &=2\ln s +\cO(\log|\log s|) \\
  &= (2+o(1)) \log s\\.
\end{align*}
We therefore find that
\[
B\sqrt\omega = (1+o(1)) \frac s {\sqrt{-2\ln s}}.
\]
Thus as $t_*\to T$, one has 
\[
\psi =\bigl(1+\cO(B^{-2}+\varepsilon)\bigr)
      (\tfrac12\sqrt{n-1}+o(1))  \frac{s}{\sqrt{-\log s}}
\]
for all $t\in[t_*,T]$ and $x_*$ for which $s=s(x_*,t_*)$ satisfies
(\ref{eq:B-curve}).  Since $B$ and $\varepsilon$ can be chosen
arbitrarily, this implies
\begin{lemma}\label{lem:outer-asymptotics}
  For any $\delta>0$ there exist $B>0$ and $ \bar t<T$
  such that

  \[
  (\tfrac12\sqrt{n-1}-\delta) \frac{s}{\sqrt{-\log s}}
  \leq\psi
  \leq (\tfrac12\sqrt{n-1}+\delta) \frac{s}{\sqrt{-\log s}}
  \]
  in the region $(T-t)|\log(T-t)|\leq s^2/B^2$, $\bar t\leq t\leq T$.
\end{lemma}

\section{Formal asymptotics for general neckpinches}
\label{NotRound}

The usual starting point for obtaining matched asymptotic expansions near a
stationary solution of a nonlinear \textsc{pde} is to linearize around that
solution.  This method encounters two obstacles when
applied to the Ricci flow.

The first challenge is the well-known fact that the Ricci flow is only weakly
parabolic.  This deficiency stems from its invariance under the full
diffeomorphism group of the underlying manifold. We overcome it by a variant of
the DeTurck trick \cite{DT1, DT2}.  (See \S \ \ref{ArbitraryPerturbations}.)
Our choice is equivalent to the Bianchi gauge, as adopted in the elliptic
context by Biquard \cite{Biq} and others.

The second challenge is more fundamental: it is the fact that the linearized
\textsc{pde} possesses a null eigenvalue, suggesting the presence of a center
manifold. This phenomenon also occurs in neck pinches of the mean curvature
flow \cite{AV}, and flat self-similar solutions to the reaction diffusion
equation $\partial_tu=\Delta u+u^p$ \cite{HV1,HV2,HV3,FK}.  We overcome it by
carrying out a higher-order variational analysis.  A heuristic description of
our method is as follows.  We will formally compute the dynamics of a center
manifold of a fixed point of a flow.  That fixed point is the \emph{cylinder
soliton }introduced in \S \ \ref{Cylinder}, and the flow is the \emph{dilated
Ricci flow }introduced in \S \ \ref{DRF}.  The center manifold will be
tangential to the null eigenspace of the linearization obtained in \S \
\ref{ArbitraryPerturbations}. A general point on the center manifold will be of
the form $\zeta+\Phi\left( \zeta\right) $, where $\zeta$ belongs to the kernel
of the linearization and $\Phi\left( \zeta\right) $ is at least quadratic in
its dependence on $\zeta$.  (See \S \ \ref{Ansatz}.)  To compute the formal
dynamics, we substitute this \emph{Ansatz }into the dilated Ricci flow, project
onto the kernel of the linearization, and compute the purely quadratic terms,
yielding an \textsc{ode} on the kernel.  The errors introduced by this method
should be smaller than the order of the solution, namely
$\mathcal{O}(\tau^{-2})$, where $\tau$ is the time scale introduced above in
Section \ref{Round}.

We pursue this method below. The computations it requires are
extensive, but the (formal) conclusion it yields is eminently simple.
It tells us that the asymptotics obtained rigorously in Section
\ref{Round} for rotationally-symmetric solutions should be stable for
fully general solutions.

\subsection{The cylinder soliton}
\label{Cylinder}

Recall that a \emph{Ricci soliton }is a tuple $\left(
\mathcal{M}^{m},g,X\right) $, where $\mathcal{M}^{m}$ is a smooth manifold, $g$
is its Riemannian metric, and $X$ is a complete vector field on
$\mathcal{M}^{m}$ such that the identity
\[
-2\operatorname*{Rc}\left(  g\right)  =\mathcal{L}_{X}g+\lambda g
\]
holds for some $\lambda\in\left\{ -1,0,1\right\} $. It is well known
that each Ricci soliton gives rise to a \emph{self-similar solution
}$\bar {g}\left( t\right) $ of the Ricci flow
\[
\frac{\partial}{\partial t}\bar{g}=-2\operatorname*{Rc}
\left(\bar{g}\right)
\]
defined by
\[
\bar{g}\left(  t\right)  =\left(  T+\lambda t\right)  \left(  \varphi
  _{t}^{\ast}g\right)  ,
\]
where $\varphi_{t}$ is the one-parameter family of diffeomorphisms
generated by the vector fields $(T+\lambda t)^{-1}X$.

Now consider the manifold $\mathbb{R}\times S^{n}$ with local coordinates
$\left( y^{0};y^{1},\ldots,y^{n}\right) $, where $y^{0}\equiv x\in \mathbb{R}$
and $\left( y^{1},\ldots,y^{n}\right) \equiv y\in S^{n}$. Let $\hat{g}\equiv
g_{\operatorname*{can}}$ denote the round metric of unit radius on $S^{n}$.
Define a product metric $g$ on $\mathbb{R}\times S^{n}$ by
\begin{equation}
  g=dx\otimes dx+2\left(  n-1\right)  \hat{g}, \label{CylinderMetric}
\end{equation}
noting that its Ricci curvature satisfies
\[
\operatorname*{Rc}\left(  g\right)  =\left(  n-1\right)  \hat{g}.
\]
Let $X$ be the vector field on $\mathbb{R}\times S^{n}$ defined by
\begin{equation}
  X=\operatorname*{grad}\left(  \frac{x^{2}}{4}\right)  =\frac{x}{2}
  \frac{\partial}{\partial x},
\end{equation}
noting that the Lie derivative of $g$ with respect to $X$ is
\[
\mathcal{L}_{X}g=dx\otimes dx.
\]
It is then easy to see that
\begin{equation}
  -2\operatorname*{Rc}\left(  g\right)  =\mathcal{L}_{X}g-g,
  \label{CylinderSoliton}
\end{equation}
hence that $\left(  \mathbb{R}\times S^{n},g,X\right)  $ is a shrinking
gradient Ricci soliton. We call this the \emph{cylinder soliton.}

For later use, we observe that the Levi-Civita connection of $g$ is
given in local coordinates $\left( y^{0};y^{1},\ldots,y^{n}\right) $
by the Christoffel symbols
\begin{equation}
  \Gamma_{ij}^{k}=\left\{
    \begin{array}
      [c]{cl}
      \hat{\Gamma}_{ij}^{k} & \text{if }1\leq i,j,k\leq n\\
      0 & \text{otherwise.}
    \end{array}
  \right.  \label{Christoffel}
\end{equation}
Here and throughout the remainder of this paper, a hat designates a geometric
quantity associated to the round unit sphere $\left( S^{n},\hat{g}\right) $.  The
components $R_{ijk\ell}=g_{\ell m}R_{ijk}^{m}$ of the Riemannian curvature of $g$ are
\begin{equation}
  R_{ijk\ell}=\left\{
    \begin{array}
      [c]{cl}
      2\left(  n-1\right)  \hat{R}_{ijk\ell}=2\left(  n-1\right)  \left(  \hat
        {g}_{i\ell}\hat{g}_{jk}-\hat{g}_{ik}\hat{g}_{j\ell}\right)  & \text{if }1\leq
      i,j,k,\ell\leq n\\
      0 & \text{otherwise,}
    \end{array}
  \right.
\end{equation}
and those of its Ricci tensor are
\begin{equation}
  R_{ij}=\left\{
    \begin{array}
      [c]{cl}
      \hat{R}_{ij}=\left(  n-1\right)  \hat{g}_{jk} & \text{if }1\leq i,j\leq n\\
      0 & \text{otherwise.}
    \end{array}
  \right.  \label{CylinderRicci}
\end{equation}

\subsection{The dilated Ricci flow\label{DRF}}

Motivated by the rigorous results obtained in Section \ref{Round}, we
want to study finite-time local singularities of the Ricci flow that
are modeled on the cylinder soliton. Accordingly, let $\left(
\mathcal{M}^{m},G\left( t\right) \right) $ be a solution of the Ricci
flow
\[
\frac{\partial}{\partial t}G=-2\operatorname*{Rc}\left(  G\right)
\]
that exists for $0\leq t<T<\infty$. Given any vector field $X$ on
$\mathcal{M}^{m}$, let $\varphi_{t}$ be the family of diffeomorphisms
solving
\[
\frac{\partial}{\partial t}\varphi\left(  \cdot,t\right)  =\frac{1}
{T-t}X\left(  \varphi\left(  \cdot,t\right)  \right)  .
\]
Then the \emph{blow-up of }$G$\emph{ at }$T$\emph{ modified by }$X$ is
the metric $g\left( t\right) $ defined by
\begin{equation}
  G\left(  t\right)  =\left(  T-t\right)  \left(  \varphi_{t}^{\ast}g\left(
      t\right)  \right)  .
\end{equation}
Observing that
\[
\frac{\partial}{\partial t}\left(  \varphi_{t}^{\ast}g\right)  =\frac{1}
{T-t}\varphi_{t}^{\ast}\left(  \mathcal{L}_{X}g\right)  +\varphi_{t}^{\ast
}\left(  \frac{\partial}{\partial t}g\right)  ,
\]
one computes that
\[
\varphi_{t}^{\ast}\left(  -2\operatorname*{Rc}\left(  g\right)  \right)
=-2\operatorname*{Rc}\left(  G\right)  =\frac{\partial}{\partial t}
G=\varphi_{t}^{\ast}\left(  -g+\mathcal{L}_{X}g+\left(  T-t\right)
  \frac{\partial}{\partial t}g\right)  .
\]
Thus if one regards $g$ as a function of the rescaled time variable
\begin{equation}
  \tau=\log\frac{1}{T-t},
\end{equation}
then $\left(  \mathcal{M}^{m},g\left(  \tau\right)  ,X\right)  $ becomes a
solution of the \emph{dilated Ricci flow:}
\begin{equation}
  \frac{\partial}{\partial\tau}g=-2\operatorname*{Rc}\left(  g\right)
  -\mathcal{L}_{X}g+g. \label{DilatedRicciFlow}
\end{equation}

The following is an immediate consequence of equations (\ref{CylinderSoliton})
and (\ref{DilatedRicciFlow}).

\begin{lemma}
  The cylinder soliton $\left(  \mathbb{R}\times S^{n},g,X\right)  $ is a
  stationary solution of the dilated Ricci flow.
\end{lemma}

\subsection{First and second order variation formulas}

If $h$ is a symmetric $\left(  2,0\right)  $-tensor on a Riemannian manifold
$\left(  \mathcal{M}^{m},g\right)  $, its Lichnerowicz Laplacian $\Delta
_{\ell}h$ is defined by
\begin{equation}
  \Delta_{\ell}h_{ij}=\Delta h_{ij}+2R_{ipqj}h^{pq}-R_{i}^{k}h_{kj}-R_{j}
  ^{k}h_{ik},
\end{equation}
where $\Delta h$ is the rough Laplacian. Its divergence $\delta h$ is defined
by
\begin{equation}
  \left(  \delta h\right)  _{j}=-\nabla^{i}h_{ij}.
\end{equation}

The following variation formulas are proved by direct calculation.

\begin{lemma}
  \label{Variation}Let $\left(  \mathcal{M}^{m},g\right)  $ be a Riemannian
  manifold; let $h$ be a symmetric $\left(  2,0\right)  $-tensor on
  $\mathcal{M}^{m}$; and let $\varepsilon$ be small enough that
  \[
  \tilde{g}=g+\varepsilon h
  \]
  is a Riemannian metric.

  \begin{enumerate}
  \item In local coordinates, the Christoffel symbols $\tilde{\Gamma}$ of the
    Levi-Civita connection of $\tilde{g}$ are
    \begin{align*}
      \tilde{\Gamma}_{ij}^{k}= 
      \Gamma_{ij}^{k}
      &  +\frac{\varepsilon}{2}
               \left(
                 \nabla_{i}h_{j}^{k}+\nabla_{j}h_{i} ^{k}-\nabla^{k}h_{ij}
               \right) \\
      &  -\frac{\varepsilon^{2}}{2}h_{\ell}^{k}\left(  \nabla_{i}h_{j}^{\ell}
        +\nabla_{j}h_{i}^{\ell}-\nabla^{\ell}h_{ij}\right)  +\mathcal{O}\left(
        \varepsilon^{3}\right)  .
    \end{align*}

  \item In local coordinates, the Riemann tensor $\widetilde{\operatorname*{Rm}
    }$ of $\tilde{g}$ is
    \begin{align*}
      \tilde{R}_{ijk}^{\ell} =R_{ijk}^{\ell}
      & +\frac{\varepsilon}{2}
      \Bigl(
        \nabla_{i}\nabla_{k}h_{j}^{\ell}
        -\nabla _{i}\nabla^{\ell}h_{jk}
        -\nabla_{j}\nabla_{k}h_{i}^{\ell}
        +\nabla_{j} \nabla^{\ell}h_{ik} \\
      & \qquad +R_{ijm}^{\ell}h_{k}^{m}-R_{ijk}^{m}h_{m}^{\ell}
      \Bigr) \\
      & +\frac{\varepsilon^{2}}{2}h_{m}^{\ell}
        \bigl(
          \nabla_{i}\nabla^{m} h_{jk}
          -\nabla_{j}\nabla^{m}h_{ik}
          -\nabla_{i}\nabla_{k}h_{j}^{m}
          +\nabla _{j}\nabla_{k}h_{i}^{m}
        \bigr) \\
      & +\frac{\varepsilon^{2}}{4}\left(
        \nabla_{i}h_{k}^{m}+\nabla_{k}h_{i}
        ^{m}-\nabla^{m}h_{ik}\right) \left(
        \nabla_{j}h_{m}^{\ell}+\nabla^{\ell
        }h_{jm}-\nabla_{m}h_{j}^{\ell}\right) \\
      & -\frac{\varepsilon^{2}}{4}\left(
        \nabla_{j}h_{k}^{m}+\nabla_{k}h_{j}
        ^{m}-\nabla^{m}h_{jk}\right) \left(
        \nabla_{i}h_{m}^{\ell}+\nabla^{\ell
        }h_{im}-\nabla_{m}h_{i}^{\ell}\right) \\
      & +\frac{\varepsilon^{2}}{2}h_{p}^{\ell}\left(
        R_{ijk}^{q}h_{q}^{p} -R_{ijq}^{p}h_{k}^{q}\right)
      +\mathcal{O}\left( \varepsilon^{3}\right) .
    \end{align*}

  \item In local coordinates, the Ricci tensor $\widetilde{\operatorname*{Rc}}$
    of $\tilde{g}$ is
    \begin{align*}
      \tilde{R}_{ij}=R_{ij}
      &  -\frac{\varepsilon}{2}\left[  \Delta_{\ell}h_{ij}+\nabla_{i}\nabla
        _{j}H+\nabla_{i}\left(  \delta h\right)  _{j}+\nabla_{j}\left(  \delta
          h\right)  _{i}\right] \\
      &  +\frac{\varepsilon^{2}}{2}h^{pq}\left(  \nabla_{p}\nabla_{q}h_{ij}
        +\nabla_{i}\nabla_{j}h_{pq}-\nabla_{i}\nabla_{p}h_{qj}-\nabla_{j}\nabla
        _{p}h_{qi}\right) \\
      &  +\frac{\varepsilon^{2}}{2}\left(  \frac{1}{2}\nabla_{i}h_{p}^{q}\nabla
        _{j}h_{q}^{p}+\nabla_{p}h_{i}^{q}\nabla^{p}h_{qj}-\nabla_{p}h_{i}^{q}
        \nabla_{q}h_{j}^{p}\right) \\
      &  +\frac{\varepsilon^{2}}{2}\left(  \frac{1}{2}\nabla_{k}H+\left(  \delta
          h\right)  _{k}\right)  \left(  \nabla_{i}h_{j}^{k}+\nabla_{j}h_{i}^{k}
        -\nabla^{k}h_{ij}\right) \\
      &  +\frac{\varepsilon^{2}}{2}h^{pq}\left(  h_{i}^{k}R_{jpkq}+h_{j}^{k}
        R_{ipkq}+2h_{q}^{k}R_{ipkj}\right)  +\mathcal{O}\left(  \varepsilon
        ^{3}\right)  .
    \end{align*}

  \item If $X$ and $Y$ are vector fields on $\mathcal{M}^{m}$ and $\tilde
    {X}=X+\varepsilon Y$, then
    \[
    \mathcal{L}_{\tilde{X}}\tilde{g}=\mathcal{L}_{X}g+\varepsilon\left(
      \mathcal{L}_{Y}g+\mathcal{L}_{X}h\right)  +\varepsilon^{2}\mathcal{L}_{Y}h;
    \]
    in particular,
    \begin{align*}
      \left(  \mathcal{L}_{\tilde{X}}\tilde{g}\right)_{ij}  =\;
      & \nabla_{i} X_{j}+\nabla_{j}X_{i}\\
      &  +\varepsilon\left[  \left(  \nabla_{i}Y_{j}+\nabla_{j}Y_{i}\right)
        +\left(  X^{k}\nabla_{k}h_{ij}+\nabla_{i}X^{k}h_{kj}+\nabla_{j}X^{k}
          h_{ik}\right)  \right] \\
      &  +\varepsilon^{2}\left(  Y^{k}\nabla_{k}h_{ij}+\nabla_{i}Y^{k}h_{kj}
        +\nabla_{j}Y^{k}h_{ik}\right)  .
    \end{align*}

  \end{enumerate}
\end{lemma}

\subsection{Arbitrary perturbations of a stationary solution}
\label{ArbitraryPerturbations}

Let $\left(\mathcal{M}^{m},g,X\right)$ be any stationary solution of the
dilated Ricci flow (\ref{DilatedRicciFlow}), and let $h$ be an arbitrary
symmetric $\left(  2,0\right)  $-tensor such that the perturbation
\begin{equation}
  \tilde{g}=g+h
\end{equation}
is a Riemannian metric on $\mathcal{M}^{m}$. Denote the trace of $h$ with
respect to $g$ by
\begin{equation}
  H=\operatorname*{tr}\!_{g}h.
\end{equation}
Let $Y$ be the vector field metrically dual to the $1$-form $\frac{1}
{2}dH+\delta h$, namely
\begin{equation}
  Y=\left(  \frac{1}{2}dH+\delta h\right)  ^{\sharp}, \label{DeTurck}
\end{equation}
and form the perturbation
\begin{equation}
  \tilde{X}=X+Y
\end{equation}
of the vector field $X$. The role of $Y$ is to implement a DeTurck
trick \cite{DT1, DT2} that makes the linearization of the dilated
Ricci flow (\ref{DilatedRicciFlow}) strictly parabolic. (In the
elliptic context, this choice is frequently called the Bianchi gauge
\cite{Biq}.) Now applying Lemma \ref{Variation} gives the following
result.

\begin{lemma}
  The metric $\tilde{g}\left( \tau\right) $ will be a solution of the
  dilated Ricci flow (\ref{DilatedRicciFlow}) for the vector field
  $\tilde{X}$ if and only if $h$ evolves by the nonlinear system
  \begin{equation}
    \frac{\partial}{\partial\tau}h=\mathcal{A}h+\mathcal{Q}\left(  h\right)
    +\mathcal{C}\left(  h\right)  , \label{PerturbedDRF}
  \end{equation}
  where $\mathcal{A}$ is the linear elliptic operator defined by
  \begin{equation}
    \mathcal{A}h=\Delta_{\ell}h-\mathcal{L}_{X}h+h, \label{FirstFormForA}
  \end{equation}
  the quadratic term $\mathcal{Q}\left( h\right) $ is defined by
  \begin{subequations}
    \label{Quadratic}
    \begin{align}
      \left[ \mathcal{Q}\left( h\right) \right] _{ij} \;=\;
      &h_{i}^{p}\nabla
      _{j}\nabla_{q}h_{p}^{q}+h_{j}^{p}\nabla_{i}\nabla_{q}h_{p}^{q}\\
      & +h^{pq}\left(
        \nabla_{i}\nabla_{p}h_{qj}+\nabla_{j}\nabla_{p}h_{qi}
        -\nabla_{p}\nabla_{q}h_{ij}-\nabla_{i}\nabla_{j}h_{pq}\right) \\
      & -\frac{1}{2}\left(
        h_{i}^{k}\nabla_{j}\nabla_{k}H+h_{j}^{k}\nabla
        _{i}\nabla_{k}H+\nabla_{i}h_{j}^{k}\nabla_{k}H+\nabla_{j}h_{i}^{k}\nabla
        _{k}H\right) \\
      &
      +\nabla_{i}h_{j}^{p}\nabla_{q}h_{p}^{q}+\nabla_{j}h_{i}^{p}\nabla_{q}
      h_{p}^{q}+\nabla_{p}h_{i}^{q}\nabla_{q}h_{j}^{p}\\
      &
      -\nabla_{p}h_{i}^{q}\nabla^{p}h_{jq}-\frac{1}{2}\nabla_{i}h_{p}^{q}
      \nabla_{j}h_{q}^{p}\\
      & +h^{pq}\left( h_{i}^{k}R_{jpqk}+h_{j}^{k}R_{ipqk}-2h_{q}^{k}
        R_{ipkj}\right) ,
    \end{align}
    and $\mathcal{C}\left( h\right) $ is at least third-order in $h$
    and its covariant derivatives.
  \end{subequations}
\end{lemma}

The cubic term $\mathcal{C}\left(  h\right)  $ will not be important in the
formal asymptotic analysis that follows.

\subsection{The linearization at the cylinder soliton}

In the case that the stationary solution of the dilated Ricci flow is the
cylinder soliton $\left(  \mathbb{R}\times S^{n},g,X\right)  $ introduced in
\S \ \ref{Cylinder}, the linear operator $\mathcal{A}$ defined in
(\ref{FirstFormForA}) can be written in a more useful form. To see this, we
again work in the coordinate system $\left(  y^{0};y^{1},\ldots,y^{n}\right)
$, with $y^{0}\equiv x\in\mathbb{R}$ and $\left(  y^{1},\ldots,y^{n}\right)
\equiv y\in S^{n}$. Recalling formula (\ref{CylinderRicci}), it is easy to see
that
\[
h_{ij}-\left(  R_{i}^{k}h_{kj}+R_{j}^{k}h_{ik}\right)  =\left\{
  \begin{array}
    [c]{cl}
    h_{00} & \text{if }i=j=0\\
    h_{i0}/2 & \text{if }1\leq i\leq n\\
    h_{0j}/2 & \text{if }1\leq j\leq n\\
    0 & \text{if }1\leq i,j\leq n.
  \end{array}
\right.
\]
This observation allows one to write the Lie derivative of $h$ with respect to
$X$ as
\begin{align*}
  \left(  \mathcal{L}_{X}h\right)  _{ij}  &  =\frac{x}{2}\frac{\partial
  }{\partial x}h_{ij}+\frac{1}{2}h_{ij}\left(  \delta_{i}^{0}+\delta_{j}
    ^{0}\right) \\
  &  =\frac{x}{2}\frac{\partial}{\partial x}h_{ij}+h_{ij}-\left(  R_{i}
    ^{k}h_{kj}+R_{j}^{k}h_{ik}\right)  .
\end{align*}
Then using the simple identity
\[
R_{i\ell}R_{jk}-R_{ik}R_{j\ell}=\frac{n-1}{2}R_{ijk\ell},
\]
one reaches the following conclusion.

\begin{lemma}
  On the cylinder soliton, the linear operator $\mathcal{A}$ defined in
  (\ref{FirstFormForA}) is given by
  \begin{equation}
    \left(  \mathcal{A}h\right)  _{ij}=\Delta h_{ij}-\frac{x}{2}\frac{\partial
    }{\partial x}h_{ij}+\frac{4}{n-1}\left(  R_{ij}R_{pq}-R_{iq}R_{pj}\right)
    h^{pq}, \label{SecondFormForA}
  \end{equation}
  where $\Delta h$ denotes the rough Laplacian.
\end{lemma}

\subsection{An upper bound for the spectrum of the linearization}
Since the manifold $\mathbb{R}\times S^{n}$ is noncompact, it is not necessarily the
case that the formally elliptic operator $\cA$ has discrete point spectrum.  For our
operator $\cA$ the first order term $-\frac12x\partial_x$ makes this happen, as we
will now show.

We consider the bundle $\cS_2$ of $(2,0)$-tensors on $\R\times S^n$ and the space
$C^\infty_0(\cS_2)$ of smooth compactly supported sections of this bundle.  The
operator $\cA$ maps this space to itself.  We define the inner product
\[
(h,k) \doteqdot 
\int_\R\int_{S^n} h_{ij}k_{ij}e^{-x^2/4}\; d\theta\;dx
\]
where $d\theta$ is the volume form on $S^n$.  We write $\|h\| =
\surd{(h,h)}$ for the corresponding norm.  The completion of
$C^\infty_0(\cS_2)$ with respect to this norm is $L^2(\cS_2;
e^{-x^2/4}\; d\theta\;dx)$, which we abbreviate to $L^2$. (See \cite{EF}, for instance.) 

\begin{lemma}
  \label{lem:A-good}
  The closure of the densely defined operator $\cA$ is a self-adjoint operator with
  compact resolvent.  $\cA$ is bounded from above.  Its spectrum consists of discrete
  finite multiplicity point spectrum.
\end{lemma}

\begin{proof}
One can write $\cA$ as
\begin{equation}
  \label{eq:A-symmetric}
  \bigl(\cA h\bigr)_{ij} = 
  e^{x^2/4} \frac\partial{\partial x}
  \left\{
    e^{-x^2/4}\frac{\partial h_{ij}}{\partial x}
  \right\}
  +\tfrac1{2(n-1)}\hat\Delta h_{ij}
  +V_{ijpq}h^{pq}
\end{equation}
in which $\hat\Delta$ is the rough Laplacian in the direction tangential to the
spheres $\{x\}\times S^n $ and
\[
V_{ijpq} \doteqdot \tfrac 4{n-1} (R_{ij}R_{pq} - R_{iq}R_{jp}).
\]
By integration by parts, one finds that for any $h,k\in
C^\infty_0(\cS_2)$ one has
\begin{align}
  \label{eq:A-quadform}
  (h,-\cA &k) =  \\ \notag
  &\int_\R\int_{S^n}\bigl\{
      \partial_x h_{ij}\partial_x k_{ij}
      +\tfrac1{2(n-1)}\hat\nabla_\ell h_{ij}\hat\nabla_\ell k_{ij}
      -V_{ijpq}h_{ij}k_{pq}
  \bigr\}
  e^{-x^2/4}\; d\theta\;dx,
\end{align}
from which it is evident that $\cA$ is symmetric.  Boundedness of $V_{ijpq}$ implies
that there is a constant $C_0$ such that
\begin{equation}
  \label{eq:A-bound}
  (h,-\cA h) \geq 
   \int_\R\int_{S^n} |\nabla h|^2 e^{-x^2/4}\; d\theta\;dx
     -C_0 \|h\|^2
\end{equation}
holds for all $h\in C^\infty_0(\cS_2)$.

Let $H^1=H^1(\cS_2)$ be the completion of $C^\infty_0(\cS_2)$ for the norm
corresponding to the inner product
\[
\langle h, k\rangle  \doteqdot 
\int_\R\int_{S^n} \bigl\{\nabla_ih_{pq}\nabla_ik_{pq} + h_{pq}k_{pq}\bigr\}
e^{-x^2/4}\; d\theta\;dx.
\]
(Again, see \cite{EF}.) Then the standard Hilbert space arguments show that for all $\lambda>C_0$, the
equation $\lambda h-\cA h = k$ has a distributional solution $h\in H^1$ for any $k\in
L^2$.  It follows that the closure of $\cA$ is a self-adjoint operator in $L^2$ whose
domain
\[
\mathrm{D}(\cA) = \{(\lambda-\cA)^{-1}k \mid k\in L^2\}
\]
is contained in $H^1$. (Here it makes no difference which $\lambda>C_0$ one chooses.)
In fact, the inequality \eqref{eq:A-bound} implies that the $H^1$ norm is equivalent
with the norm $||| h |||^2 = \langle h, (\lambda-\cA)h\rangle$, so that $H^1$ is the
form domain of $\cA$, i.e.~$H^1=\mathrm{D}(\sqrt{\lambda-\cA})$

The inequality (\ref{eq:A-bound}) tells us that the spectrum of $\cA$ is
contained in the interval $(-\infty, C_0]$.

To see that $\cA$ has compact resolvent, we recall Lemma~\ref{lem:sigmaW-bound},
which implies
\[
\| xh \|^2 \leq C\int_\R\int_{S^n} 
  \{|\partial_xh|^2 + |h|^2\} e^{-x^2/4}\;d\theta\;dx
  \leq C\|h\|_{H^1}^2.
\]
Hence
\[
\|h\|^2+\|xh\|^2+\|\nabla h\|^2 \leq C (h, (\lambda-\cA) h).
\]
This inequality together with the Rellich-Kondrachov theorem imply that the imbedding
$H^1\subset L^2$ is compact.  Since $H^1$ is the form domain of $\cA$, we conclude
that $\cA$ indeed has compact resolvent, so that its spectrum is pure point spectrum.

\end{proof}

\textsc{Remark. } More generally, the linearization of the modified Ricci flow operator at any
gradient Ricci soliton $(\mathcal{M}^{n},g)$ will be self-adjoint in
$L^{2}(\mathcal{M}^{n},g;\,e^{-\varphi}\,d\mu)$ if $d\mu$ is the measure
naturally associated to $g$ and one takes $\varphi$ to be the soliton
potential function.

\subsection{A decomposition into invariant subspaces}

In order to compute the spectrum of $\cA$, we split the space $L^2 (\mathcal{S}_2)$
into a number of $\cA$-invariant subspaces.

It is a standard fact that any $\left( 2,0\right) $-tensor $h$ on a Riemannian
manifold $\left( \mathcal{M}^{m},g\right) $ admits the decomposition
\[
h=h^{\vee}+h^{\circ},
\]
where
\[
h^{\vee}=\left(  \operatorname*{tr}\!_{g}h\right)  \,g
\]
is a multiple of the metric, and
\[
h^{\circ}=h-\frac{1}{m}h^{\vee}
\]
is trace free.

One can exploit the product structure of the cylinder soliton $\left(
  \mathbb{R}\times S^{n},g,X\right) $ to obtain a more refined decomposition.  Again
using coordinates $\left( y^{0};y^{1},\ldots,y^{n}\right) \ $ and adopting the
convention that Roman indices lie in the range $0\dots n$ while Greek indices lie in
$1\dots n$, we decompose a given symmetric $\left( 2,0\right) $-tensor $h$ on
$\mathbb{R}\times S^{n}$ as
\begin{subequations}
  \label{CanonicalDecomposition}
  \begin{align}
    h  &  =M\,dx\otimes dx+N_{\alpha}\left(  dx\otimes dy^{\alpha}+dy^{\alpha
      }\otimes dx\right)  +\left(  O_{\alpha\beta}^{\vee}+O_{\alpha\beta}^{\circ
      }\right)  \,dy^{\alpha}\otimes dy^{\beta}\\
    &  =M\,dx\otimes dx+N_{\alpha}\left(  dx\otimes dy^{\alpha}+dy^{\alpha}\otimes
      dx\right)  +\left(  P\hat{g}_{\alpha\beta}+Q_{\alpha\beta}\right)
    \,dy^{\alpha}\otimes dy^{\beta}.
  \end{align}
  Here one should regard $M$ as a scalar function on $\mathbb{R}\times S^{n}$,
  $N\left(  x,\cdot\right)  $ as a one-parameter family of $1$-forms defined on
  the spheres $\left\{  x\right\}  \times S^{n}$, $P$ as the scalar function
\end{subequations}
\[
P=\frac{2\left(  n-1\right)  }{n}\left(  H-M\right)  ,
\]
and $Q\left(  x,\cdot\right)  $ as the one-parameter family of trace-free
$\left(  2,0\right)  $-tensors defined on the spheres $\left\{  x\right\}
\times S^{n}$ by
\[
Q_{\alpha\beta}=h_{\alpha\beta}-P\hat{g}_{\alpha\beta}.
\]

The utility of this decomposition is that it exhibits a set of invariant
subspaces for $\mathcal{A}$. To demonstrate this, it will be helpful to fix
additional notation. Let $\hat{\Delta}$ denote the rough Laplacian of the
canonical sphere metric $\hat{g}$. Let $-\hat{\Delta}_{dR}$ denote its
Hodge--de Rham Laplacian, recalling that $\hat{\Delta}_{dR}$ acts on a
$1$-form $N$ by
\[
\hat{\Delta}_{dR}N_{\alpha}=-\left[  \left(  d\delta+\delta d\right)
  N\right]  _{\alpha}=\hat{\Delta}N_{\alpha}-\hat{R}_{\alpha}^{\beta}N_{\beta
}=\hat{\Delta}N_{\alpha}-\left(  n-1\right)  N_{\alpha}.
\]
Let $\hat{\Delta}_{\ell}$ denote the Lichnerowicz Laplacian of $\hat{g}$,
which acts on a trace-free tensor $Q$ by
\[
\hat{\Delta}_{\ell}Q_{\alpha\beta}=\hat{\Delta}Q_{\alpha\beta}-2nQ_{\alpha
  \beta}.
\]
Finally, let $\mathcal{B}$ denote the differential operator defined by
\begin{equation} \label{DefineHermiteOperator}
  \mathcal{B}=\frac{\partial^{2}}{\partial x^{2}}-\frac{x}{2}\frac{\partial
  }{\partial x}
  = e^{x^2/4}\frac\partial{\partial x} e^{-x^2/4}\frac\partial{\partial x}.
\end{equation}
Using these conventions, one computes using equation (\ref{SecondFormForA})
that the decomposition of $\mathcal{A}h$ corresponding to
(\ref{CanonicalDecomposition}) is
\begin{align*}
  \left(  \mathcal{A}H\right)  _{00}  &  =\left[  \mathcal{B}+\frac{1}{2\left(
        n-1\right)  }\hat{\Delta}\right]  M\\
  \left(  \mathcal{A}H\right)  _{\alpha0}  &  =\left[  \mathcal{B}+\frac
    {1}{2\left(  n-1\right)  }\hat{\Delta}\right]  N_{\alpha}=\left[
    \mathcal{B}+\frac{1}{2\left(  n-1\right)  }\hat{\Delta}_{dR}+\frac{1}
    {2}\right]  N_{\alpha}\\
  \left(  \mathcal{A}H\right)  _{\alpha\beta}^{\vee}  &  =\left\{  \left[
      \mathcal{B}+\frac{1}{2\left(  n-1\right)  }\hat{\Delta}+1\right]  P\right\}
  \hat{g}_{\alpha\beta}\\
  \left(  \mathcal{A}H\right)  _{\alpha\beta}^{\circ}  &  =\left[
    \mathcal{B}+\frac{1}{2\left(  n-1\right)  }\hat{\Delta}-\frac{1}{n-1}\right]
  Q_{\alpha\beta}.
\end{align*}
This calculation proves the following result.

\begin{lemma}
  \label{SpectralDecomposition}On the cylinder soliton, the operator
  $\mathcal{A}$ given by (\ref{SecondFormForA}) may be represented schematically
  with respect to the decomposition (\ref{CanonicalDecomposition}) as
  \begin{equation}
    \mathcal{A}:
    \begin{pmatrix}
      M\\
      \mathstrut\\
      N\\
      \mathstrut\\
      P\\
      \mathstrut\\
      Q
    \end{pmatrix}
    \mapsto
    \begin{pmatrix}
      \left[  \mathcal{B}+\frac{1}{2\left(  n-1\right)  }\hat{\Delta}\right]  M\\
      \mathstrut\\
      \left[  \mathcal{B}+\frac{1}{2\left(  n-1\right)  }\hat{\Delta}_{dR}+\frac
        {1}{2}\right]  N\\
      \mathstrut\\
      \left[  \mathcal{B}+\frac{1}{2\left(  n-1\right)  }\hat{\Delta}+1\right]  P\\
      \mathstrut\\
      \left[  \mathcal{B}+\frac{1}{2\left(  n-1\right)  }\hat{\Delta}-\frac{1}
        {n-1}\right]  Q
    \end{pmatrix}
    . \label{diagonalizeA}
  \end{equation}
  The only contribution to the spectrum of $\mathcal{A}$ comes from its four
  component operators: $\mathcal{B}+\frac{1}{2\left(  n-1\right)  }\hat{\Delta}$
  and $\mathcal{B}+\frac{1}{2\left(  n-1\right)  }\hat{\Delta}+1$ acting on
  scalar functions, $\mathcal{B}+\frac{1}{2\left(  n-1\right)  }\hat{\Delta
  }_{dR}+\frac{1}{2}$ acting on $1$-forms, and $\mathcal{B}+\frac{1}{2\left(
      n-1\right)  }\hat{\Delta}-\frac{1}{n-1}$ acting on trace-free $\left(
    2,0\right)  $-tensors.
\end{lemma}

\subsection{The spectrum of the linearization\label{Spectrum}}

The decomposition obtained in Lemma \ref{SpectralDecomposition} allows us to
analyze the spectrum of $\mathcal{A}$ using separation of variables. To
illustrate the idea, suppose that
\[
M\left(  x,y\right)  =h\left(  x\right)  F\left(  y\right)  ,
\]
where $h$ is an eigenfunction of $\mathcal{B}$ with eigenvalue $\lambda$, and
$F$ is an eigenfunction of $\frac{1}{2\left(  n-1\right)  }\hat{\Delta}$ with
eigenvalue $\mu$. Then
\[
\left[  \mathcal{B}+\frac{1}{2\left(  n-1\right)  }\hat{\Delta}\right]
M=\left(  \lambda+\mu\right)  M.
\]
This line of argument leads easily to the following observation.

\begin{lemma}
  \label{SeparateVariables}All eigenvalues of
  $\mathcal{A}:\mathrm{D}(\cA)\to L^{2}$ have the form $\lambda+\mu$, where
  $\lambda$ is an eigenvalue of $\mathcal{B}$ acting on $\mathfrak{h}
  =L^{2}(\mathbb{R};e^{-x^{2}/4})$, and $\mu$ is an eigenvalue of either
  $\frac{1}{2\left( n-1\right) }\hat{\Delta}$ or $\frac{1}{2\left( n-1\right)
  }\hat{\Delta}+1$ acting on $C^{\infty}\left( S^{n}\right) $, of $\frac{1}{2\left(
      n-1\right) }\hat{\Delta}_{dR}+\frac{1}{2}$ acting on $\Omega_{1}\left(
    S^{n}\right) $, or of $\frac{1}{2\left( n-1\right) } \hat{\Delta}-\frac{1}{n-1}$
  acting on the space of smooth trace-free $\left( 2,0\right) $-tensors on $S^{n}$.
\end{lemma}

To make this result useful, we need some simple facts about the component
operators that appear in (\ref{diagonalizeA}).

\medskip

\textsc{The operator }$\mathcal{B}=\frac{\partial^{2}}{\partial x^{2}}
-\frac{x}{2}\frac{\partial}{\partial x}$\textsc{ acting on scalar functions.
}The Hermite polynomials $\left\{  h_{k}\right\}  _{k=0}^{\infty}$ constitute
a complete orthogonal family for the operator $\mathcal{B}$ on the weighted
Hilbert space
\[
\mathfrak{h}=L^{2}\left(  \mathbb{R},e^{-x^{2}/4}\,dx\right)  .
\]
Moreover, one has
\[
\left(  \mathcal{B}+1\right)  h_{k}=\left(  1-\frac{k}{2}\right)  h_{k},
\]
so that the spectrum of $\mathcal{B}$ is $\left\{  0,-\frac{1}{2},-1,-\frac
  {3}{2},\ldots\right\}  $. (Notice that $\mathfrak{h}$ is a larger space than
we defined in Section \ref{Round}, since we not restricting to even functions here.)

\textsc{The operator }$\frac{1}{2\left(  n-1\right)  }\hat{\Delta}$\textsc{
  acting on scalar functions. }It is well known that the spectrum of
$\hat{\Delta}$ acting on $C^{\infty}\left(  S^{n}\right)  $ is $\left\{
  \lambda_{k}\right\}  _{k\geq0}$, where
\[
\lambda_{k}=-k\left(  n+k-1\right)  .
\]
In particular, the only nonnegative eigenvalue of $\frac{1}{2\left(
    n-1\right)  }\hat{\Delta}$ is $\lambda_{0}=0$. Its eigenspace consists of the
constant functions.

\textsc{The operator }$\frac{1}{2\left(  n-1\right)  }\hat{\Delta}+1$\textsc{
  acting on scalar functions. }If $k\geq2$, then
\[
\lambda_{k}+2\left(  n-1\right)  =\left(  2-k\right)  \left(  n-1\right)
-k^{2}\leq-k^{2}<0.
\]
Hence the only possible nonnegative eigenvalues of $\frac{1}{2\left(
    n-1\right)  }\hat{\Delta}+1$ are
\[
\frac{\lambda_{0}+2\left(  n-1\right)  }{2\left(  n-1\right)  }=1
\]
and
\[
\frac{\lambda_{1}+2\left(  n-1\right)  }{2\left(  n-1\right)  }=\frac{1}
{2}\frac{n-2}{n-1}\in\left[  0,\frac{1}{2}\right)  .
\]
The eigenspace corresponding to $1$ consists of constants. The eigenspace
corresponding to $\frac{1}{2}\frac{n-2}{n-1}$ consists of the spherical
harmonics: the restrictions to $S^{n}$ of the linear functions on
$\mathbb{R}^{n+1}$.

\textsc{The operator }$\frac{1}{2\left(  n-1\right)  }\hat{\Delta}_{dR}
+\frac{1}{2}$\textsc{ acting on }$1$\textsc{-forms. }It is shown in \cite{B}
that the spectrum of $\hat{\Delta}_{dR}$ acting on the space $\Omega
_{1}\left(  S^{n}\right)  $ of smooth $1$-forms is
\[
\left\{  \lambda_{k}\right\}  _{k\geq1}\bigcup\left\{  \lambda_{k}
  +2-n\right\}  _{k\geq1}.
\]
If $k\geq1$, then $\lambda_{k}+n-1\leq-1$. Because $\lambda_{k}+2-n\leq
\lambda_{k}$ in all dimensions $n\geq2$, it follows that $\frac{1}{2\left(
    n-1\right)  }\hat{\Delta}_{dR}+\frac{1}{2}$ is negative definite.

\textsc{The operator }$\frac{1}{2\left(  n-1\right)  }\hat{\Delta}-\frac
{1}{n-1}$\textsc{ acting on trace-free }$\left(  2,0\right)  $
\textsc{-tensors. }The operator $\frac{1}{2\left(  n-1\right)  }\hat{\Delta
}-\frac{1}{n-1}$ is clearly negative definite.

Combining these observations with Lemma \ref{SeparateVariables} leads to the
following conclusion.

\begin{lemma}
  \label{SpectralSummary}On the cylinder soliton, the only possible non-negative
  eigenvalues of the linearization $\mathcal{A}$ are $0$, $\frac{1}{2}\frac
  {n-2}{n-1}$, $\frac{1}{2}$, and $1$.

  \begin{itemize}
  \item The eigenspace corresponding to $1$ consists of constant multiples of
    $\hat{g}$.

  \item The eigenspace corresponding to $\frac{1}{2}$ consists of multiples of
    $\hat{g}$ that are linear in $x\in\mathbb{R}$.

  \item The eigenspace corresponding to $\frac{1}{2}\frac{n-2}{n-1}$ consists of
    spherical harmonics.

  \item The eigenspace corresponding to $0$ consists of constant functions and
    of multiples of $\hat{g}$ that are quadratic in $x\in\mathbb{R}$. (When $n=2$,
    the spherical harmonics also belong to this eigenspace.)
  \end{itemize}
\end{lemma}

\medskip

\textsc{Remark. }The results above can be obtained in another way. Since all
geometric data of the metric (\ref{CylinderMetric}) are independent of
$x\in\mathbb{R}$, one has the simple commutator
\[
\left[  \mathcal{A},\frac{\partial}{\partial x}\right]  =\frac{1}{2}
\frac{\partial}{\partial x}.
\]
So if $\mathcal{A}h=\lambda h$, a straightforward induction argument shows
that
\[
\mathcal{A}\left(  \frac{\partial^{k}}{\partial x^{k}}h\right)  =\left(
  \lambda+\frac{k}{2}\right)  \left(  \frac{\partial^{k}}{\partial x^{k}
  }h\right)
\]
for all integers $k\geq0$. Since Lemma \ref{lem:A-good} proves that the
spectrum of $\mathcal{A}$ is bounded from above, it follows readily that for
every eigentensor $h$ of $\mathcal{A}$, there exist some integer
$\kappa=\kappa\left(  h\right)  \geq0$ and a family of $\left(  2,0\right)
$-tensors $\hat{h}_{1},\ldots,\hat{h}_{\kappa}$ independent of $x\in
\mathbb{R}$ such that
\[
h=\sum_{k=0}^{\kappa}x^{k}\hat{h}_{k}.
\]

\subsection{Our \emph{Ansatz}\label{Ansatz}}

The analysis in \S \ \ref{Spectrum} shows that any eigenvalues of
$\mathcal{A}$ corresponding to trace-free (i.e.~non-rotationally symmetric)
perturbations of the cylinder soliton are all strictly negative. As we
explained in the introduction to this section, the only perturbations relevant
to our formal center-manifold computation are those corresponding to null
eigenvalues of the linearization. Therefore, we adopt the \emph{Ansatz }that
\[
h=M\,dx\otimes dx+N_{\alpha}\left(  dx\otimes dy^{\alpha}+dy^{\alpha}\otimes
  dx\right)  +\left(  P\hat{g}_{\alpha\beta}+Q_{\alpha\beta}\right)
\,dy^{\alpha}\otimes dy^{\beta}
\]
has the special form
\begin{equation}
  h\left(  x,y,\tau\right)  =u\left(  x,\tau\right)  \,dx\otimes dx+\left\{
    v\left(  x,\tau\right)  +w\left(  y,\tau\right)  \right\}  \,\hat
  {g},\label{h-Ansatz}
\end{equation}
where once again, $\tau$ denotes rescaled time, $x$ denotes a point in
$\mathbb{R}$, and $y$ denotes a point in $S^{n}$. Specifically --- continuing
to let $\left\{  h_{k}\right\}  _{k=0}^{\infty}$ denote the Hermite
polynomials --- we assume that $M\equiv u$, where
\begin{equation}
  u\left(  x,\tau\right)  =u_{0}\left(  \tau\right)  \,h_{0}\left(  x\right)
  =u_{0}\left(  \tau\right)  .\label{IA1}
\end{equation}
We further assume that $N\equiv0$, that $Q\equiv0$, and that $P\equiv v+w$,
where
\begin{equation}
  v\left(  x,\tau\right)  =v_{2}\left(  \tau\right)  \,h_{2}\left(  x\right)
  =v_{2}\left(  \tau\right)  \,\left(  x^{2}-2\right)  \label{IA2}
\end{equation}
and
\begin{equation}
  w\left(  y,\tau\right)  =\left\{
    \begin{array}
      [c]{cl}
      \sum_{i=1}^{3}\omega_{i}\left(  \tau\right)  \,\Omega_{i}\left(  y\right)   &
      \text{if }n=2\\
      0 & \text{otherwise.}
    \end{array}
  \right.  \label{IA3}
\end{equation}
Here $\left(  \Omega_{1},\Omega_{2},\Omega_{3}\right)  $ is a basis for the
space of spherical harmonics on $S^{2}$. (Recall that $w$ corresponds to a
positive eigenvalue of the linearization whenever the \emph{total dimension
}$n+1$ is at least $4$.)

To implement our formal center-manifold analysis, we shall study the flow
\begin{equation}
  \frac{\partial}{\partial\tau}h=\mathcal{A}h+\mathcal{Q}\left(  h\right)
  \label{QRF}
\end{equation}
that models equation (\ref{PerturbedDRF}) up to second order. To do so
requires us to analyze the quadratic term $\mathcal{Q}\left(  h\right)  $ at
the cylinder soliton in a manner analogous to what was done above for the
linearization $\mathcal{A}$. To simplify the notation, we again use
coordinates $\left(  y^{0};y^{1},\ldots,y^{n}\right)  $ and assume that Roman
indices lie in the range $0\dots n$ while Greek indices lie in $1\dots n$. If
$f$ is a smooth function, we further adopt the convention that $f_{x}=\partial
f/\partial x$ and $f_{\alpha}=\partial f/\partial y^{\alpha}$.

We shall begin with computations that are more general than what we need at
the moment. To wit, we assume only that $h$ has the form given by
(\ref{h-Ansatz}), without imposing the specific assumptions (\ref{IA1})--(\ref{IA3}).
One then verifies readily that all first covariant derivatives
of the tensor $h$ vanish except
\begin{align*}
  \nabla_{0}h_{00}  &  =u_{x}\\
  \nabla_{0}h_{\alpha\beta}  &  =v_{x}\,\hat{g}_{\alpha\beta}\\
  \nabla_{\gamma}h_{\alpha\beta}  &  =w_{\gamma}\,\hat{g}_{\alpha\beta}.
\end{align*}
All second covariant derivatives of $h$ vanish except
\begin{align*}
  \nabla_{0}\nabla_{0}h_{00}  &  =u_{xx}\\
  \nabla_{0}\nabla_{0}h_{\alpha\beta}  &  =v_{xx}\,\hat{g}_{\alpha\beta}\\
  \nabla_{\gamma}\nabla_{\delta}h_{\alpha\beta}  &  =\hat{\nabla}_{\gamma}
  \hat{\nabla}_{\delta}w\,\hat{g}_{\alpha\beta}.
\end{align*}
Moreover, the only nonvanishing derivatives of the trace
\[
H=g^{ij}h_{ij}=u+\frac{n}{2\left(  n-1\right)  }\left(  v+w\right)
\]
are
\begin{align*}
  \nabla_{0}H  &  =u_{x}+\frac{n}{2\left(  n-1\right)  }v_{x}\\
  \nabla_{\gamma}H  &  =\frac{n}{2\left(  n-1\right)  }w_{\gamma}
\end{align*}
and
\begin{align*}
  \nabla_{0}\nabla_{0}H  &  =u_{xx}+\frac{n}{2\left(  n-1\right)  }v_{xx}\\
  \nabla_{\gamma}\nabla_{\delta}H  &  =\frac{n}{2\left(  n-1\right)  }
  \hat{\nabla}_{\gamma}\hat{\nabla}_{\delta}w.
\end{align*}

Using these formulas, a straightforward calculation reveals that the first
component of the quadratic term
\begin{align*}
  \left[  \mathcal{Q}\left(  h\right)  \right]  _{00}    =
  &2h_{00}\nabla _{0}\nabla_{0}h_{00}-h^{\alpha\beta}\nabla_{0}\nabla_{0}h_{\alpha\beta}
  -h_{00}\nabla_{0}\nabla_{0}H\\
  &  -\nabla_{0}h_{00}\nabla_{0}H
     +\frac{3}{2}\left(  \nabla_{0}h_{00}\right)^{2}
     -\frac{1}{2}\nabla_{0}h_{\alpha}^{\beta}\nabla_{0}h_{\beta}^{\alpha}
\end{align*}
may be written as
\begin{subequations}
  \label{Q00}
  \begin{align}
    \left[  \mathcal{Q}\left(  h\right)  \right]  _{00}    =
    &u{}u_{xx}-\frac{n}{4\left(n-1\right)^{2}}v_{xx}\left[2(n-1)u+v+w\right] \\
    &  +\frac{1}{2}u_{x}^{2}
       -\frac{n}{8(n-1)^{2}}v_{x}  \left[4(n-1)u_{x} + v_{x}\right]  .
  \end{align}
\end{subequations}
Then using the fact that
\[
g^{\gamma\delta}\hat{\nabla}_{\gamma}\hat{\nabla}_{\delta}w=\frac{1}{2\left(
    n-1\right)  }\hat{\Delta}w=-\frac{n}{2\left(  n-1\right)  }w,
\]
one computes that the remaining components
\begin{align*}
  \left[  \mathcal{Q}\left(  h\right)  \right]  _{\alpha\beta}    =&h_{\alpha
  }^{\gamma}\nabla_{\beta}\nabla_{\delta}h_{\gamma}^{\delta}+h_{\beta}^{\gamma
  }\nabla_{\alpha}\nabla_{\delta}h_{\gamma}^{\delta}-h^{00}\nabla_{0}\nabla
  _{0}h_{\alpha\beta}\\
  &  +h^{\gamma\delta}\left(  \nabla_{\alpha}\nabla_{\gamma}h_{\delta\beta
    }+\nabla_{\beta}\nabla_{\gamma}h_{\delta\alpha}-\nabla_{\gamma}\nabla_{\delta
    }h_{\alpha\beta}-\nabla_{\alpha}\nabla_{\beta}h_{\gamma\delta}\right) \\
  &  -\frac{1}{2}\left(  h_{\alpha}^{\gamma}\nabla_{\beta}\nabla_{\gamma
    }H+h_{\beta}^{\gamma}\nabla_{\alpha}\nabla_{\gamma}H+\nabla_{\alpha}h_{\beta
    }^{\gamma}\nabla_{\gamma}H+\nabla_{\beta}h_{\alpha}^{\gamma}\nabla_{\gamma
    }H\right) \\
  &  +\nabla_{\alpha}h_{\beta}^{\gamma}\nabla_{\delta}h_{\gamma}^{\delta}
  +\nabla_{\beta}h_{\alpha}^{\gamma}\nabla_{\delta}h_{\gamma}^{\delta}
  +\nabla_{\gamma}h_{\alpha}^{\delta}\nabla_{\delta}h_{\beta}^{\gamma}\\
  &  -\nabla_{0}h_{\alpha}^{\gamma}\nabla^{0}h_{\beta\gamma}-\nabla_{\delta
  }h_{\alpha}^{\gamma}\nabla^{\delta}h_{\beta\gamma}-\frac{1}{2}\nabla_{\alpha
  }h_{\gamma}^{\delta}\nabla_{\beta}h_{\delta}^{\gamma}\\
  &  +h^{\gamma\delta}\left(  h_{\alpha}^{\eta}R_{\beta\gamma\delta\eta
    }+h_{\beta}^{\eta}R_{\alpha\gamma\delta\eta}-2h_{\delta}^{\eta}R_{\alpha
      \gamma\eta\beta}\right)
\end{align*}
are given by
\begin{subequations}
  \label{Qab}
  \begin{align}
    \left[  \mathcal{Q}\left(  h\right)  \right]  _{\alpha\beta}    =
    &-\frac {1}{2\left(  n-1\right)  }\left\{  2\left(  n-1\right)  uv_{xx}+v_{x}
      ^{2}\right\}  \,\hat{g}_{\alpha\beta}\\
    &  -\frac{1}{4\left(  n-1\right)  ^{2}}\left\{  \left\vert \hat{\nabla
        }w\right\vert _{\hat{g}}^{2}-nw\left(  v+w\right)  \right\}  \,\hat{g}
    _{\alpha\beta}\\
    &  -\frac{n-2}{8\left(  n-1\right)  ^{2}}\left\{  4\left(  v+w\right)
      \hat{\nabla}_{\alpha}\hat{\nabla}_{\beta}w+3w_{\alpha}w_{\beta}\right\}  .
  \end{align}
\end{subequations}

\subsection{The inner layer $\left\vert x\right\vert =o\left(  \sqrt{\tau
    }\right)  $}

We now derive a formal profile of a solution
\[
\tilde{g}\left(  x,y,\tau\right)  =g\left(  x,y\right)  +h\left(
  x,y,\tau\right)
\]
near a singularity modeled on the cylinder soliton $g$. We assume that
(\ref{h-Ansatz}) and (\ref{IA1})--(\ref{IA3}) hold. These assumptions force
$h$ to belong to the kernel of the linearization of equation
(\ref{PerturbedDRF}), because
\[
\mathcal{A}h=\left[  \mathcal{B}u\right]  \,dx\otimes dx+\left[  \left(
    \mathcal{B}+1\right)  v\right]  \,\hat{g}=0.
\]
Recall that $w$ vanishes by assumption if $n>2$. In this context, (\ref{Q00})
becomes
\begin{equation}
  \left[  \mathcal{Q}\left(  h\right)  \right]  _{00}=-\left\{  v_{2}w+\frac
    {n}{\left(  n-1\right)  ^{2}}\left[  \left(  n-1\right)  u_{0}v_{2}+v_{2}
      ^{2}\right]  \right\}  h_{0}-\left\{  \frac{n}{\left(  n-1\right)  ^{2}}
    v_{2}^{2}\right\}  h_{2}. \label{Q00-refined}
\end{equation}
And (\ref{Qab}) reduces to
\[
\left[  \mathcal{Q}\left(  h\right)  \right]  _{\alpha\beta}=\mathcal{P}
\,\hat{g}_{\alpha\beta},
\]
where
\begin{equation}
  \mathcal{P}=\left\{  \frac{1}{2}w^{2}-\frac{1}{4}\left\vert \hat{\nabla
      }w\right\vert _{\hat{g}}^{2}-2u_{0}v_{2}-\frac{4}{n-1}v_{2}^{2}\right\}
  h_{0}+\left\{  \frac{1}{2}v_{2}w-\frac{2}{n-1}v_{2}^{2}\right\}  h_{2}.
  \label{Qab-refined}
\end{equation}

In order efficiently to derive \textsc{ode} for the functions $u_{0}$, $v_{2}
$, and $\omega_{1},\omega_{2},\omega_{3}$ from equations (\ref{Q00-refined})
and (\ref{Qab-refined}), it is helpful to adopt some additional notation: if
$f$ is a function of $y\in S^{n}$, we define
\[
\left\Vert f\right\Vert _{S^{n}}^{2}=\int_{S^{n}}f\left(  y\right)
^{2}\,d\hat{\mu};
\]
and if $X$ and $Y$ are symmetric $\left(  2,0\right)  $-tensors$\ $on
$\mathbb{R}\times S^{n}$, we define
\[
\llangle X,Y\rrangle\,=\int_{-\infty}^{\infty}\left\langle X,Y\right\rangle
_{x}e^{-x^{2}/4}\,dx,
\]
where
$\left\langle X,Y \right\rangle _{x}=\int_{\{x\}\times S^{n}}\left\langle X,Y\right\rangle
_{\hat{g}}\,d\hat{\mu}$. By using the fact that $\left\{  h_{k}\right\}
_{k=0}^{\infty}$ is an orthogonal family for $\mathfrak{h}$, it is then easy
to see that
\begin{align*}
  u_{0}^{\prime}\left(  \tau\right)  \cdot\left\Vert h_{0}\right\Vert
  _{\mathfrak{h}}^{2}\left\Vert 1\right\Vert _{S^{n}}^{2}  &  =\frac{d}{dt}\llangle
  h,h_{0}\,dx^{2}\rrangle\\
  &  =\,\llangle\mathcal{Q}\left(  h\right)  ,h_{0}\,dx^{2}\rrangle\\
  &  =-\frac{n}{\left(  n-1\right)  ^{2}}\left[  \left(  n-1\right)  u_{0}
    v_{2}+v_{2}^{2}\right]  \cdot\left\Vert h_{0}\right\Vert _{\mathfrak{h}}
  ^{2}\left\Vert 1\right\Vert _{S^{n}}^{2}
\end{align*}
and
\begin{align*}
  v_{2}^{\prime}\left(  \tau\right)  \cdot n\left\Vert h_{2}\right\Vert
  _{\mathfrak{h}}^{2}\left\Vert 1\right\Vert _{S^{n}}^{2}  &  =\frac{d}{dt}\llangle
  h,h_{2}\,\hat{g}\rrangle\\
  &  =\,\llangle\mathcal{Q}\left(  h\right)  ,h_{2}\,\hat{g}\rrangle\\
  &  =-\frac{2}{n-1}v_{2}^{2}\cdot n\left\Vert h_{2}\right\Vert _{\mathfrak{h}
  }^{2}\left\Vert 1\right\Vert _{S^{n}}^{2}.
\end{align*}
Similarly, for any $i\in\left\{  1,2,3\right\}  $, we find when $n=2$ that
\[
\omega_{i}^{\prime}\left(  \tau\right)  \cdot2\left\Vert h_{0}\right\Vert
_{\mathfrak{h}}^{2}\left\Vert \Omega_{i}\right\Vert _{S^{2}}^{2}=\frac{d}
{dt}\llangle h,\Omega_{i}\,\hat{g}\rrangle\,=\,\llangle\mathcal{Q}\left(  h\right)
,\Omega_{i}\,\hat{g}\rrangle\,=0.
\]
This calculations let us derive a formal profile of $h$ in the inner layer
$\left\vert x\right\vert =o\left(  \sqrt{\tau}\right)  $.

\begin{lemma}
  \label{InnerLayerLemma}For perturbations of the form
  \[
  h=\left\{  u_{0}\left(  \tau\right)  \,h_{0}\left(  x\right)  \right\}
  \,dx\otimes dx+\left\{  v_{2}\left(  \tau\right)  \,h_{2}\left(  x\right)
    +\sum_{i=1}^{3}\omega_{i}\left(  \tau\right)  \,\Omega_{i}\left(  y\right)
  \right\}  \,\hat{g},
  \]
  the flow (\ref{QRF}) is equivalent to the system of \textsc{ode}
  \label{CenterManifold}
  \begin{align}
    \frac{d}{d\tau}u_{0}  &  =-\frac{n}{\left(  n-1\right)  ^{2}}\left[  \left(
        n-1\right)  u_{0}v_{2}+v_{2}^{2}\right] \label{CenterManifold1}\\
    \frac{d}{d\tau}v_{2}  &  =-\frac{2}{n-1}v_{2}^{2}\label{CenterManifold2}\\
    \frac{d}{dt}\omega_{i}  &  =0\qquad\qquad\left(  i=1,2,3\right)
    \label{CenterManifold3}
  \end{align}
  whose solutions, up to an error term of $o(\tau^{-1})$ as $\tau\rightarrow
  \infty$, are \label{InnerLayer}
  \begin{align}
    u_{0}\left(  \tau\right)   &  =\left\{
      \begin{array}
        [c]{cl}
        -\frac{\log\left(  2\tau\right)  }{2\tau}+\frac{c}{\tau} & \text{if }n=2\\
        \mathstrut & \mathstrut\\
        -\frac{n}{2\left(  n-2\right)  }\frac{1}{\tau} & \text{otherwise}
      \end{array}
    \right. \\
    v_{2}\left(  \tau\right)   &  =\frac{n-1}{2\tau}\\
    \omega_{i}\left(  \tau\right)   &  =\omega_{i}\left(  0\right)  \qquad
    \qquad\qquad\qquad\left(  i=1,2,3\right)  .
  \end{align}
  Formally, this is the profile of $h$ in the region $\left\vert x\right\vert
  =o\left(  \sqrt{\tau}\right)  $.
\end{lemma}

\begin{proof}
  If $v_{2}$ is not identically zero, equation (\ref{CenterManifold2}) has the
  explicit solution
  \[
  v_{2}\left(  \tau\right)  =\frac{n-1}{c+2\tau}.
  \]
  Then equation (\ref{CenterManifold1}) becomes
  \[
  \frac{d}{d\tau}u_{0}\left(  \tau\right)  +\left(  \frac{n}{c+2\tau}\right)
  u_{0}=-\frac{n}{\left(  c+2\tau\right)  ^{2}}.
  \]
  An integrating factor for this linear \textsc{ode} is $\left(  c+2\tau\right)
  ^{n/2}$, whence one obtains
  \[
  u_{0}\left(  \tau\right)  =\left\{
    \begin{array}
      [c]{cl}
      -\frac{\log\left(  c+2\tau\right)  }{c+2\tau}+\frac{C}{c+2\tau} & \text{if
      }n=2\\
      \mathstrut & \mathstrut\\
      -\frac{n}{n-2}\frac{1}{c+2\tau}+\frac{C}{\left(  c+2\tau\right)  ^{n/2}} &
      \text{otherwise.}
    \end{array}
  \right.
  \]

\end{proof}

\medskip

\textsc{Remark. }To compare the formal results above with those obtained
rigorously in Section \ref{Round}, observe that our choices of dilating
factors imply that the quantity $a_{2}\left(  \tau\right)  $ appearing in
Section \ref{Round} and the quantity $v_{2}\left(  \tau\right)  $ appearing
above are related by
\[
a_{2}\left(  \tau\right)  \approx\frac{v_{2}\left(  \tau\right)  }{4\left(
    n-1\right)  }\approx\frac{1}{8\tau},
\]
because in the inner layer, we have $x\approx\sigma$ and
\begin{align*}
  &  \left\{  2\left(  n-1\right)  +v_{2}\left(  \tau\right)  h_{2}\left(
      x\right)  +o(\tau^{-1})\right\}  \,\hat{g}\\
  &  \qquad\qquad\qquad\qquad=2\left(  n-1\right)  \left[  1+a_{2}\left(
      \tau\right)  h_{2}\left(  \sigma\right)  +o(\tau^{-1})\right]  ^{2}\hat{g}\\
  &  \qquad\qquad\qquad\qquad=\left\{  2\left(  n-1\right)  +4\left(
      n-1\right)  a_{2}\left(  \tau\right)  h_{2}\left(  \sigma\right)  +o(\tau
    ^{-1})\right\}  \,\hat{g}.
\end{align*}

\textsc{Remark. }By modifying $\tilde{g}$ by an initial conformal
diffeomorphism, one may make all $\omega_{i}$ vanish. Therefore, throughout
the remainder of this paper, we will assume that $\omega_{i}\equiv0$ for
$i=1,2,3$.

\subsection{The intermediate layer $\left\vert x\right\vert =\mathcal{O}
  \left(  \sqrt{\tau}\right)  $}

In order to study the region where $\left\vert x\right\vert =\mathcal{O}
\left(  \sqrt{\tau}\right)  $, we replace the coordinates $\left(
  x,\tau\right)  $ with $\left(  \xi,\tau\right)  $, where
\begin{equation}
  \xi=\frac{x}{\sqrt{\tau}}. \label{define-xi}
\end{equation}
If $f$ is a smooth function of $\left(  x,\tau\right)  $, we define $F\left(
  \xi,\tau\right)  =f\left(  x,\tau\right) .$ Then one has the formulas
\[
f_{\tau}=F_{\tau}-\frac{\xi}{2\tau}F_{\xi}\qquad\text{and}\qquad f_{x}
=\frac{1}{\sqrt{\tau}}F_{\xi}.
\]

Let
\[
U\left(  \xi,\tau\right)  =u\left(  x,\tau\right)  \qquad\text{and}\qquad
V\left(  \xi,\tau\right)  =v\left(  x,\tau\right)  ,
\]
where $u$ and $v$ are the quantities that appear in the simple \emph{Ansatz
}(\ref{h-Ansatz}). Then substituting formulas (\ref{Q00}) and (\ref{Qab}) into
the \textsc{pde} (\ref{QRF}), one computes that
\[
U_{\tau}=-\frac{\xi}{2}U_{\xi}+\frac{1}{\tau}\left\{
  \begin{array}
    [c]{c}
    \left(  1+U\right)  U_{\xi\xi}+\frac{\xi}{2}U_{\xi}+\frac{1}{2}U_{\xi}^{2}\\
    \mathstrut\\
    -\frac{n}{4\left(  n-1\right)  ^{2}}V_{\xi\xi}\left[  2\left(  n-1\right)
      U+V\right] \\
    \mathstrut\\
    -\frac{n}{8\left(  n-1\right)  ^{2}}V_{\xi}\left[  4\left(  n-1\right)
      U_{\xi}+V_{\xi}\right]
  \end{array}
\right\}
\]
and
\[
V_{\tau}=V-\frac{\xi}{2}V_{\xi}+\frac{1}{\tau}\left\{  V_{\xi\xi}\left(
    1-U\right)  +\frac{\xi}{2}V_{\xi}-\frac{1}{2\left(  n-1\right)  }V_{\xi}
  ^{2}\right\}  .
\]
Recalling Lemma \ref{InnerLayerLemma}, we will construct an approximate
solution by considering the equations
\[
U_{\tau}=-\frac{\xi}{2}U_{\xi}+o(\tau^{-1})
\]
and
\[
V_{\tau}=V-\frac{\xi}{2}V_{\xi}+\mathcal{O}(\tau^{-2}).
\]
These equations suggest that the solution is modeled by the time-independent
profiles $U\left(  \xi\right)  \equiv C_{1}$ and $V\left(  \xi\right)
=C_{2}\xi^{2}$, respectively. Matching this layer with the results of Lemma
\ref{InnerLayerLemma} determines $C_{1}$ and $C_{2}$, yielding the following result.

\begin{lemma}
  \label{MiddleLayerLemma}In the intermediate layer $\left\vert x/\sqrt{\tau
    }\right\vert =\left\vert \xi\right\vert =\mathcal{O}\left(  1\right)  $, the
  formal solution is modeled by
  \[
  u\left(  x,\tau\right)  \approx U\left(  \frac{x}{\sqrt{\tau}}\right)
  \equiv0\qquad\text{and}\qquad v\left(  x,\tau\right)  \approx V\left(
    \frac{x}{\sqrt{\tau}}\right)  =\frac{n-1}{2}\frac{x^{2}}{\tau}.
  \]
  To wit,
  \[
  \tilde{g}\approx dx\otimes dx+2\left(  n-1\right)  \left[  1+\left(
      \xi/2\right)  ^{2}\right]  \,\hat{g}.
  \]

\end{lemma}

\subsection{The outer layer $\left\vert x\right\vert =o\left(  1/
    {\sqrt{T-t}}\right)  $}

To describe the outer layer, we first introduce the \textquotedblleft
blown-down\textquotedblright\ coordinate
\begin{equation}
  s=e^{-\tau/2}x=\sqrt{T-t}\cdot x,
\end{equation}
which should be compared to the metric distance $s$ defined in Section
\ref{Round}.

If one were to follow the method above by defining $\bar{U}\left(
  s,\tau\right)  =u\left(  x,\tau\right)  $ and $\bar{V}\left(  s,\tau\right)
=v\left(  x,\tau\right)  $, an easy computation would show that
\[
\bar{U}_{\tau}=\mathcal{O}\left(  T-t\right)  \qquad\text{and}\qquad\bar
{V}_{\tau}=V+\mathcal{O}\left(  T-t\right)  .
\]
Unfortunately, the steady-state solutions to these equations fail to convey
enough useful information.

Instead, we proceed as follows. Fix a large number $\Xi>0$. For each $s$, let
$t_{s}$ be the time such that $\left\vert \xi\right\vert =\Xi$ at $t_{s}$,
where $\xi$ is defined by (\ref{define-xi}). Let $\tau_{s}=\log\frac
{1}{1-t_{s}}$. Then
\[
\Xi=\left\vert \xi\right\vert =\frac{\left\vert x\right\vert }{\sqrt{\tau_{s}
  }}=\frac{\left\vert s\right\vert }{\sqrt{\left(  T-t_{s}\right)  \tau_{s}}}.
\]
This shows that
\[
\tau_{s}=\log\left(  \Xi^{2}\right)  +\log\frac{1}{s^{2}}+\log\tau_{s},
\]
hence that
\begin{equation}
  \tau_{s}=\left[  2+o\left(  1\right)  \right]  \log\frac{1}{\left\vert
      s\right\vert }\qquad\text{as}\qquad\left\vert s\right\vert \rightarrow
  0.\label{tau-s}
\end{equation}
Now we \textquotedblleft blow down\textquotedblright\ the solution of the
dilated (\textquotedblleft blown-up\textquotedblright) Ricci flow by defining
\[
\underline{u}\left(  s,\tau\right)  =e^{-\tau}u\left(  x,\tau\right)
\qquad\text{and}\qquad\underline{v}\left(  s,\tau\right)  =e^{-\tau}v\left(
  x,\tau\right)  .
\]
Then at time $t_{s}$, we compute that
\[
\underline{v}\left(  s,\tau_{s}\right)  =e^{-\tau_{s}}V\left(  \frac{x}
  {\sqrt{\tau_{s}}}\right)  =\frac{n-1}{2}\frac{s^{2}}{\tau_{s}}.
\]
When $\left\vert s\right\vert $ is small and $\Xi$ is large, the undilated
solution $\underline{v}$ cannot change much in the short time $t_{s}\leq t<T$.
Hence (\ref{tau-s}) implies that as $\left\vert s\right\vert \rightarrow0$,
one has
\begin{align*}
  \underline{v}\left(  s,\tau_{s}\right)   &  =\frac{n-1}{2}\frac{s^{2}}{\left[
      2+o(1)\right]  \log\frac{1}{\left\vert s\right\vert }}\left[  1+o(1)\right]
  \\
  &  =\frac{n-1}{2}\frac{s^{2}}{\log\frac{1}{s^{2}}}\left[  1+o(1)\right]  .
\end{align*}
Since $\underline{u}\left(  s,\tau\right)  =e^{-\tau}U\left(  x/\sqrt{\tau
  }\right)  \equiv0$, we have obtained the following result.

\begin{lemma}
  In the outer layer $\left\vert x\right\vert =o\left(  \frac{1}{\sqrt{T-t}
    }\right)  $, the formal solution is modeled by
  \[
  u\left(  x,\tau\right)  \approx\frac{1}{T-t}\,\underline{u}\left(
    s,\tau\right)  \qquad\text{and}\qquad v\left(  x,\tau\right)  \approx\frac
  {1}{T-t}\,\underline{v}\left(  s,\tau\right)  ,
  \]
  where
  \[
  \underline{u}\left(  s,\tau\right)  \equiv0\qquad\text{and}\qquad\underline
  {v}\left(  s,\tau\right)  =\frac{n-1}{2}\frac{s^{2}}{\log\frac{1}{s^{2}}
  }\left[  1+o(1)\right]  .
  \]
  To wit,
  \[
  v\left(  x,\tau\right)  \approx\frac{n-1}{2}\frac{x^{2}}{\tau+\log\frac
    {1}{x^{2}}}.
  \]

\end{lemma}

\end{document}